\documentclass[dvips,preprint,11pt]{imsart}
\usepackage{color}
\usepackage{ amsmath, amsfonts, amssymb, amsthm, amscd}
\usepackage[T1]{fontenc}
\usepackage[english]{babel}
\usepackage{palatino}
\usepackage[overload]{empheq}
\usepackage{latexsym,dsfont,textcomp,amsmath}
\usepackage[english]{babel}
\usepackage[latin1]{inputenc}
\usepackage{enumerate,indentfirst}
\usepackage[colorlinks,citecolor=blue,urlcolor=blue]{hyperref}
\usepackage{amssymb}
\usepackage{amsthm}
\usepackage{hyperref}
\usepackage{graphicx}

\usepackage[latin1]{inputenc}

\def\gg{{\mathbb G}}

\def\P{\mathbb P}
\def\Es{\mathbb E}

\def\Q{\mathbb Q}
\def \os \hat{o}

\def \sd {S^{\delta} }
\def \yi {y^i} 
\def \yun {y^{1}}
\def \yd {y^d}
\def \y {y }
\def \yiun {y^{i,1}}
\def \yideux {y^{i,2}}
\def \yk {y^{k}}
\def \lqi {\ell ^i}
\def \vh {v}
\def \vi {v^i} 

 \def \lqsi {\ell ^{*,i} }
\def\Q{\mathbb Q}
\def \yud {Y^{12}} 
\def \uudt {{  U}^{12}}
\def \uudb {{\bar  U}^{12}}
\def \finproof {\hfill $  \Box$}

\def\cro#1{\langle #1\rangle}
\def \wt {\widetilde}
\usepackage[T1]{fontenc}
\usepackage[english]{babel}

\newtheorem{theoreme}{Theorem}
\newtheorem{lemma}{Lemma}
\newtheorem{Hypothese}{Assumption A}
\newtheorem{proposition}{Proposition}

\newtheorem{definition}{Definition}

\newtheorem{remarque}{Remark}
\newtheorem{exemple}{Example}

\newcommand{\cG}{\ensuremath{\mathcal G}}

\newcommand{\N}{\mathbb{N}}

\newcommand{\be}{\begin{equation}}
\newcommand{\ee}{\end{equation}}
\newcommand{\beq}{\begin{eqnarray}}
\newcommand{\eeq}{\end{eqnarray}}
\newcommand{\beqa}{\begin{eqnarray*}}
\newcommand{\eeqa}{\end{eqnarray*}}

\newcommand{\R}{\mathbb{R}}

\newcommand{\ind}{{\bf {1}}}

\newcommand{\ced}{\end{proof}}
\newcommand{\mb}[1]{\makebox{#1}}

 \def\mb{\textcolor{blue} }

\begin{document}
\begin{frontmatter}
\title{Robust utility maximization problem in model with jumps  and unbounded claim}
\date{}
\runtitle{}

\author{\fnms{Monique}
 \snm{JEANBLANC}\corref{}\ead[label=e1]{monique.jeanblanc@univ-evry.fr}}
\thankstext{T1}{The work of the first author is supported by the Chair \textit{risque de cr\'edit},
 F\'ed\'eration Bancaire Fran\c{c}aise}
\address{ Universit\'e d'\'Evry-Val-d'Essonne \\
Laboratoire de Math\'ematiques et Mod\'elisation d'\'Evry (LaMME)\\   UMR CNRS 8071 \\
\printead{e1}
}

\author{\fnms{Anis}
 \snm{MATOUSSI}\corref{}\ead[label=e2]{anis.matoussi@univ-lemans.fr}}
 \thankstext{t3}{Research partly supported by the Chair {\it Financial Risks} of the {\it Risk Foundation} sponsored by Soci\'et\'e G\'en\'erale, the Chair {\it Derivatives of the Future} sponsored by the {F\'ed\'eration Bancaire Fran\c{c}aise}, and the Chair {\it Finance and Sustainable Development} sponsored by EDF and Calyon }
\address{
 Université  du Mans \\
Insititut du Risque et de l'Assurance du Mans \\
Laboratoire Manceau de Math\'ematiques\\
\printead{e2}
 }

\author{\fnms{Armand}
 \snm{NGOUPEYOU}\corref{}\ead[label=e3]{ngoupeyou@beac.int}}
\address{ 
  Banque Des Etats de l'Afrique Centrale\\
  Service des Etudes\\
  \printead{e3}
}

\runauthor{M. Jeanblanc , A. Matoussi, A. Ngoupeyou}

\vspace{3mm}

\begin{abstract}
We study a problem  of utility maximization under model uncertainty  with information including  jumps.  We prove first that  the  value process of  the robust  stochastic control problem is described by the solution of a quadratic-exponential backward stochastic differential equation with jumps.  Then, we  establish a dynamic maximum principle for the optimal control of the maximization problem. The characterization of the optimal model and the optimal control (consumption-investment)  is given  via a forward-backward system which  generalizes the  result of Duffie and Skiadas \cite{DS94} and El  Karoui et al. \cite{ElkPQ01} in the case of  maximization of recursive utilities including model with jumps.
\end{abstract}

\noindent {\bf Key words~:}  Robust maximization problem,  preferences, model uncertainty, stochastic control, recursive utility, stochastic differential utility, backward stochastic differential equations, forward-backward system,  maximum principle, jump model.

\end{frontmatter}


\section{Introduction}
The utility maximization is a basic problem in mathematical finance which was introduced
by Merton \cite{mer}. Using stochastic control methods, he has exhibited  a closed
form formula for the value function and the optimal portfolio/consumption when the dynamics of the 
risky  asset  follows a geometric Brownian motion and the utility   function is  of
CRRA type. There exists  a   huge  literature on this problem  based on  two
approaches: the Bellman approach and the martingale one.
Karatzas et al. \cite{kar1}  have studied   a consumption-investment problem {in a more general case than  \cite{mer}},  taking
  into consideration  the inherent non-negativity of consumption and
bankruptcy constraint. When the risky assets are modeled by  geometric Brownian motions, the value function is determined explicitly, as well as the optimal consumption and the investment strategy, by solving the Bellman Equation and using a
verification theorem. Martingale methods were introduced by Karatzas et al. \cite{karle} and Cox and
Huang \cite{cox89}, who characterized  the optimal consumption- portfolio
policies when there are non-negativity constraints on both consumption and final
wealth. Unlike the nonlinear Bellman equation, they  gave a verification theorem which involves a linear partial
differential equation. In all these works, the authors  have assumed  that the underlying model is exactly known.
\par\medskip
Our paper deals with the problem of utility maximization from a
terminal value and an intermediate control under model
uncertainty. Uncertainty refers to the case in which a decision
maker does not know the probability distribution governing
the stochastic nature of the problem she/he is facing.  This
uncertainty is captured by using capacities or sets of probability
measures over the space of state of the world. The set of such
probability  measures  on  some measurable space $ (\Omega,
\mathcal F)$ is called by economics \emph{objectively rational
beliefs}, and each element of such set  is called beliefs on
$\mathcal F$  that the decision maker is able to justify on the
basis of the available information. The incompletness of
information is then captured by the fact that  one considers  a
 set of probability measures not reduced to a   singleton.
\par\medskip
In the mathematical finance literature,  there are two approaches
to solve robust utility maximization problems. The first
one relies on duality methods which are presented in  Quenez
\cite{Q04}, Gundel \cite{G05}, Shied and Wu \cite{SW05} and
Shied \cite{S05}. The second approach, which is  the one followed
in this paper, is based on the penalization method and the
minimization is taken over all possible models  as in  Anderson,
Hansen and Sargent  \cite{AHS03} and Hansen et al.
\cite{HSTW06}. These authors have introduced and discussed the
robust utility maximization problem when the model
uncertainty $ \Q$ is penalized by a relative entropy term with
respect to a given reference probability measure $\P$. Both papers
are cast in a Markovian setting and use mainly formal
manipulations of Hamilton-Jacobi-Bellman (HJB) equations to
provide insights about the optimal investment behaviour in these
situations.  In \cite{S03},  Skiadas  follows the same point of
view and gives the dynamics of the control problem via BSDE
in the Markovian context. More precisely, Skiadas  \cite{S03}
points out that the BSDE coincides with  the one describing a
stochastic differential utility; hence, working with a standard
expected utility under (a particular form of) model uncertainty is
equivalent to working with a corresponding stochastic differential
utility under a fixed model (see also Duffie and Epstein
\cite{DE92} and Duffie and Skiadas \cite{DS94} for more about
stochastic differential utilities).   We have to mention the
interesting works of Maccheroni,  Marinacci and Rustichini
\cite{MMR06-1, MMR06-2} who have studied preferences and
robustness by using  variational  technics.
\par\medskip
More recently, Bordigoni, Matoussi and Schweizer \cite{BMS07}  have studied  this robust problem in more general setting (non Markovian approach) by using  stochastic control technics. They  have  considered the robust maximization problem  :
 {
 \begin{equation}
 \label{pb}
 \sup_{ {\psi},c}\inf_{\Q \in {\cal Q}
 }
\textbf{U}\left(\left(\psi,c\right),\Q\right)
\end{equation}
where $  \psi$ runs through a set of random variables, $ c$
through a set of   processes and $\Q$ through a set of models
(measures), and where the  criteria    $
\textbf{U}\left(\left(\psi,c\right),\Q\right)$ is the sum of a
$\Q$-expected utility and  \textit{a penalization term} associated
with  the relative entropy. They have solved only the minimization
problem and   proved  the existence of a unique $\Q^{*}$ optimal
minimizer model. Moreover, in the case of \textit{continuous
filtration}, they  have used the dynamic programming Bellman
principle to show that the value function of the stochastic
control problem is the unique solution of  a generalized BSDE with
quadratic driver. \par\medskip Bordigoni \cite{BMS07th} has
studied partially the maximization problem by using classical
optimization arguments and assuming some conjectures. She has
derived the G\^ateaux differential of the Lagrangian associated to
the optimization problem. She has obtained necessary and
sufficient conditions that must be fulfilled for the existence
of an optimal strategy in a complete market in the case of
consumption/investment problem.\par\medskip
In  \cite{FMM11} and \cite{MMM15}, the second author have  studied the maximization part of problem \eqref{pb} by using the BSDE approach as in Duffie and Skiadas \cite{DS94} and  El Karoui, Peng and Quenez \cite{ElkPQ01} in the case of continuous filtration. We mention also that there is another approach based on Ekeland Variational principle to obtain a dynamic maximum principle for recursive utility optimization problem (see Ji and Zhou \cite{jizhou1,jizhou2}). \par\medskip
In our paper, we prove first  that in the case of \textit{discontinuous filtration} (information  including  jumps), the value process $V$ of   the stochastic minimization problem in \eqref{pb} is described by a class of quadratic-exponential  BSDE with jumps (QBSDEJs in short).
 Moreover,  we characterize the minimal optimal probability measure by means of the martingale part solution  appearing in these  QBSDEJs and we prove existence and uniqueness of  solution  of  this class of equations,  by using the related  stochastic control technics. We stress that for a given unbounded  terminal condition, the study of QBSDEs is a difficult problem, see for instance Briand and Hu \cite{BH05, BH08} and Barrieu and El Karoui \cite{BElk13}  in the continuous framework and we emphasize that adding jumps in \cite{BMS07} involves  significant difficulties in solving the related BSDEs .
Then, in order to tackle the maximization problem, we prove a
comparison theorem for this class of QBSDEJs with unbounded final
condition which allows us to prove a dynamic maximum principle for
the stochastic control problem in terms of the minimal
optimal probability measure (see Proposition 8 and Theorem 5).
These results may be considered as a generalization  of the
maximum principle proven in  \cite[Theorem 4.2]{ElkPQ01} .
Moreover, characterization of the optimal control $ (\psi^*, c^*,
\Q^*)$ as the solution of a  forward-backward system is given.
Finally, we prove in the case of logarithmic utility function  of
the control process $c$ that the optimal control $ (0, c^*, \Q^*)$
is characterized via  the unique solution of some
forward-backward stochastic differential equation.  This part of
our work  is a generalization of the dynamic maximum principle
obtained by El Karoui, Quenez and Peng \cite{ElkPQ01}, in a
framework including jumps and under model uncertainty. Our results
may also be considered as a generalization of the works of
\cite{DS94, FMM11, SS99, SS08}.

 Finally, we
have to emphasize that some classes of QBSDEJs were studied
by Becherer \cite{B06} and Morlais \cite{M09a, M09b} in the case
of bounded terminal condition (see also \cite{ABL09} ). More
recently, using the forward approach introduced in \cite{BElk13},
El Karoui, Matoussi and Ngoupeyou \cite{ElkMN11} have  obtained an
existence result of a general class of QBSDEJs with unbounded
final condition.\par\medskip

 The paper is organized as follows. In Section $2$, we describe  the
model and the stochastic control problem. In Section $3$, we characterize
the optimal model measure for the minimization problem via a solution of QBSDEJs. We prove a comparison and  a
concavity  result   for the solution of our QBSDEJ  with respect to control parameters.
 In Section $4$, we  derive the necessary and sufficient conditions that must satisfy the optimal control and then we establish the dynamic maximum principle which characterizes implicitly the optimal strategy $ (\psi^*, c^*, \Q^*)$ as a solution of a forward-backward system.
For a specific choice of utility functions, the value function  is  given in Section 5 in terms of
the optimal plan. Finally, Section 6 contains a technical proof concerning  a regularity result of QBSDEJs.

\section{The Model and the  Stochastic  Control Problem }
We consider  a filtered probability space $(\Omega,{\cal G},
\mathbb G, \P)$. All the processes   are    $\mathbb G$-adapted,
and  defined  on the time interval $[0,T]$ where $T$ is
the finite horizon.    We recall that any  special $\gg$-semimartingale $Y$ admits
a canonical decomposition $Y=Y_0+A +M^{Y,c}+M^{Y,d}$ where $A $ is
a predictable finite variation process, $M^{Y,c}$ is a continuous
martingale  and $M^{Y,d}$ is a discontinuous martingale.

\begin{Hypothese} \label{hyp} We make the following
assumptions:
\begin{enumerate}
\item For each $i=1,\dots, d$, $H^i$ is a counting
process and there exists a   positive adapted process
$\lambda^i$, called the $\P$ intensity of $H^i$,  such that  the
process  $N^i$ with
\begin{equation}
\label{martingaledefaut}
N^i_t: = H^i_t-\int_0^t\lambda^i_s ds
\end{equation} is a
martingale. We assume that the processes   $H^i, i=1,\dots, d$  have
no common jumps.
\item Any discontinuous martingale $M^{ d}$ admits a representation of the
 form $dM^{d}_t=\sum_{i=1}^d \yi _tdN^i_t $ where $\yi, i=1,\dots ,d$
are   predictable processes.
\end{enumerate}
\end{Hypothese}
\noindent
 This hypothesis is satisfied in the case where the filtration is
generated by a continuous martingale  and an inhomogeneous $d$-dimensional Poisson
process, or in the case where the counting processes are stopped after the first jump, as it is done in credit risk, as  in the following example:

\begin{exemple}\label{exemple1}
  For each $i=1,\dots, d$,  let $H^i_t=\ind_{\{\tau_i\leq
t\}}$, where $\tau_i$ is a  random time and $\P(\tau_i=\tau_j)=0, i\neq j$.  Let  $\mathbb G$ be the smallest right-continuous filtration which contains the
filtration $\mathbb F^B$  generated by a $p$ dimensional Brownian motion $B$ and the filtration generated by the
processes $H^i$. Under
   the assumption that $\P( \tau_i \in d\theta_i; i=1,\dots, d \vert
\mathcal F_t^B) \sim  \P( \tau_i \in d\theta_i; i=1,\dots, d)$, then  condition (\ref{martingaledefaut}) is satisfied and    any  local $\gg$-martingale $\zeta={(\zeta_t)}_{t\ge 0}$ admits the following decomposition: $\P$-a.s,
\begin{equation}
\label{representation} \zeta_t=\zeta_0+\int_0^t  Z_s \cdot dW_s+
\int_0^t U_s  \cdot dN_s  \qquad \forall \,  t\ge 0
\end{equation}
where $W$ is the
martingale part of the $\gg$-semimartingale $B$ (see \cite{cam}) ,   $Z := (Z^1,\dots, Z^p)$ and $ U:=(U^1,\dots,U^d)$ are
$\gg$-predictable processes. Furthermore, if $\zeta$ is square
integrable
$$\Es^{\P} \big[\int_0^T |Z_s|^2ds\big]<\infty, \quad
\sum_{i=1}^d\Es^{\P}\big[\int_0^T  |U^i_s|^2 \lambda^i_sds
\big]<\infty .$$
\end{exemple}
\noindent We denote by $|X|$   the  Euclidean norm of a vector or a row vector $X$. We give now  some notations and definitions :
 \begin{definition}{~}
 \\\\
 \noindent  $L^{\exp} $ is the space of all ${\cal{ G}}_T$-measurable random variables $X $ with
 $\Es^{\P}\left[\exp\left(\gamma \vert X  \vert\right)\right]<\infty, \quad \forall
 \gamma>0.$
 \\\\\noindent 
 $D_0^{\exp}$
 is the space of  progressively measurable processes $X={(X_t)}_{t\le T}$ with
 $$\Es^{\P}\left[\exp\left(\gamma {~ \rm{ess}\sup}_{0\le t\le T}|X_t|\right)\right]<\infty, \quad
 \forall \gamma>0.$$
 $D^{\exp}_1$ is the space of  progressively measurable processes $X={(X_t)}_{t\le T}$ such that
$$\Es^{\P}\left[\exp\left(\gamma\int_0^T |X_s|ds\right)\right]<\infty, \quad \forall \gamma >0.$$
${\cal M}_0^p $ is the space of  $\P$-martingales $M={(M_t)}_{t\le
T}$ with $M_0=0$  and $\Es^{\P}\big[\sup_{0\le t\le
T}|M_t|^p\big]<\infty. $\\\\
$\mathcal{L}^2(\lambda )$ is the space of $\R^d$-valued
predictable processes $X$ such that  $ \displaystyle \sum_{i=1}^d
\Es^{\P}\Big[\int_0^T (X^i_s)^2
\lambda ^i_sds\Big]<\infty.$\\\\
$\mathcal{H}^2  $  is the space of $\R$-valued predictable
processes $X$ such that  $ \displaystyle \Es^{\P}\Big[\int_0^T |X_s|^2
 ds\Big]<\infty.$ \\\\
$\mathcal{S}^2 $ is the space of all $\R$-valued predictable
processes $X$ such that  $ \displaystyle \Es^{\P}\Big[\sup_{ 0\leq t \leq T} |X_s|^2\Big]<\infty$.\\\\
${\cal M}_{0,loc}^c $ is the set of continuous $\P$-local
martingales.\\\\
In case it is necessary, we shall indicate the probability measure we are working with, e.g., ${\cal M}_{0,loc}^c (\Q) $ for the set of $\Q$-local martingales.
\end{definition}

\begin{definition}{~}\\\\
For any probability measure $\Q$ on $(\Omega,{\cal G}_T)$,
$$
H(\Q\vert \P):=\left\{
\begin{array}{l}
\Es^{\Q}\left[\ln \frac{d\Q}{d\P}\right] \hspace{0.3cm} \hbox{ if } \Q<\!\!<\P \hbox{ on } {\cal G}_T\\
+\infty \hspace{2 cm} otherwise
\end{array}
\right.
$$
 is  the relative entropy of $\Q$ with respect to $\P$. We
denote by ${\cal Q}_f$ (resp. ${\cal Q}_f^e$)  the space of all probability {measures}
$\Q$ on $(\Omega,{\cal G}_T)$ with $\Q<\!\!<\P$ (resp. equivalent  to $\P$) on ${\cal G}_T$
and $H(\Q\vert \P)<+\infty$.  Note that the reference probability
measure $ \mathbb{P}$ belongs to  ${\cal Q}_f^e$.
\end{definition}

\subsection{The robust optimization problem}
We define a discounting process $\sd_t:=e^{-\int_0^t
\delta_s ds}$ for all $t \in [0,T]$ where $\delta$ is a
non-negative  adapted  process. For $\Q \in {\cal Q}_f$, we denote
by $Z^{\Q}={(Z^{\Q}_t)}_{0\le t\le T}$ (a c\`adl\`ag
$\P$-martingale) its  Radon-Nikodym density with respect to $\P$.  Let $U$ be a given process (the cost process) and $\bar U_T$ a
given random variable (the terminal target). The robust utility
maximization problem ${\cal {P}}(U, \bar U_T, \beta  )$ is
 to find the infimum of $\Gamma (\Q) $ over the set ${\cal Q}_f$ where
   \begin{equation}\label{expectation2}
   \begin{split}
{\Gamma(\Q)}&:=\Es^{\Q}\left[ \int_0^T \sd _s {U}_sds+  \sd _T
\bar U_T \right]+ \beta \Es^{\Q}\left[\int_0^T \delta_s \sd _s \ln
Z^{\Q}_s ds+\sd _T\ln Z^{\Q}_T\right]\,\\
&  = :\Es^{\Q}\big[ {\mathcal U}_{0,T}^{\delta} \big] + \beta \Es^{\Q}\big[{\mathcal R}_{0,T}^{\delta} (\Q) \big]
\end{split}
\end{equation}
 The first term in the right-hand side of \eqref{expectation2} will
be linked, in the following section, to the $\Q$-expected
discounted utility from target and cost process.  The second term
is a discounted relative entropy term and $ \beta >0 $ is a  given
positive constant  which determines the strength of this penalty
term. Note that the optimal probability $\Q $ for the problem
${\cal {P}}(U, \bar U_T,\beta)$  is optimal for the minimization
problem ${\cal {P}}(U^\beta, \bar U_T^\beta,1)$ where $U^\beta =
U/\beta,\, \bar U^\beta_T=\bar U_T/\beta$, therefore, we shall
restrict our attention to the problem ${\cal {P}}(U, \bar
U_T):={\cal {P}}(U, \bar
U_T,1)$.
\begin{Hypothese} \label{hyp2} For  a more precise formulation of our problem, we make   the   following further assumptions:\\[0.2cm]
i) the discount rate $\delta$ is a positive bounded process,
more precisely  there exist   two constants $\epsilon >0$ and $ c
> 0$ such that for any $t\geq 0$,
$ 0 < \epsilon \leq \delta_t  \leq  \|\delta\|_{\infty} \leq c, $ a.s. \\[0.2cm]
ii) the cost process $U$  belongs to  $ D^{\exp}_1$ and the terminal target $\bar U_T$ is in $L^{\exp} $.\\[0.2cm]
iii) the   process $  \Lambda_t^i := \int_0^t \lambda_s^i ds $  is
  uniformly bounded, i.e.,  $ \Lambda_T^i \leq C, \,
a.s.$.
\end{Hypothese}
\begin{remarque}
The assumption $U$  belongs to  $ D^{\exp}_1$  implies that  $ \mathbb{E}^{\Q} \big[ \int_0^T |U_s| ds \big] < \infty  $ for all $\Q \in {\cal Q}_f$. Indeed, applying the following estimate: $$ xy \leq y \ln y - y + e^x, \qquad \mbox{for all} \, x \in \mathbb{R}, \; y \geq 0,$$
we get $$ \mathbb{E}^\Q \big[ \int_0^T |U_s| ds \big] = \mathbb{E}^\P \big[ Z^{\Q}_T \int_0^T |U_s| ds \big]
\leq H(\Q\vert \P) - \mathbb{E}^\P \big[ Z^{\Q}_T \big] + \mathbb{E}^\P \big[ \exp \big( \int_0^T |U_s|ds \big) \, \Big].$$
\end{remarque}
\begin{remarque}
The assumption  iii)  is a technical hypothesis needed  only  in
the proof of Theorem \ref{difdym}.
\end{remarque}
We recall the existence result of the optimal
probability measure for the minimization problem $\inf_{\Q \in
{\cal {Q}}_f } \Gamma(\Q)$ which was given in Theorem 9 and Theorem 12 in Bordigoni, Matoussi et Schweizer \cite{BMS07}:
\begin{proposition}\label{optimal}
Under Assumptions   A\ref{hyp}-A\ref{hyp2},
  there exists a unique $\Q^*$ which minimizes
$ \Gamma(\Q)$ over all $\Q\in {\cal Q}_f$:
\begin{equation}\label{objectif}\Gamma (\Q^*)=\inf_{\Q \in {\cal
{Q}}_f}  \Gamma(\Q)\end{equation}  Furthermore,  $\Q^*$ is
equivalent to $\P$,  i.e.,  $\Q^* \in {\mathcal Q}_f^e $.
\end{proposition}

 \section{The Optimal Model Measure and BSDE}
\subsection{A BSDE  description of the value process}
\label{BSDEvalue}
We use stochastic control technics to describe  the dynamics of the value process
$V$ associated with our robust optimization problem, via BSDEs.  In a  markovian framework with continuous filtration (see Skiadas \cite{S03}) or in a  continuous semimartingale setting  (see Bordigoni, Matoussi and Schweizer \cite{BMS07}), the authors  have established  that $V$ is the unique solution of a backward stochastic differential equation (BSDE) with a  quadratic driver.
 In our paper,  the BSDE associated with our control problem  (in a framework including jumps)  will  contain   quadratic and exponential terms and will be of the following form:
\begin{definition}
 A triple of processes $(Y,M^{Y,c}, \y)$ such that $Y$ is a $\P$-semimartingale,
$M^{Y,c}$ is a locally square-integrable continuous local
$\P$-martingale null at 0 and $\y=(\yun,\dots,\yd)$ an
$\mathbb{R}^d$-valued predictable locally bounded process, is called solution of   BSDEJ,  if it  satisfies:
\begin{equation}\label{eqback}
\left\lbrace
\begin{split}
&dY_t=\Big[ \sum_{i=1}^d  g( \yi_t)  \lambda^i_t- U_t+\delta_t
Y_t\Big]dt + \frac{1}{2} d{\langle M^{Y,c} \rangle}_t
+dM^{Y,c}_t+\sum_{i=1}^d { \yi _t}dN^i_t\\
&Y_T= \bar{U}_T
\end{split}
\right.
\end{equation}
where $g$ is the convex function  $g(x)=e^{-x}+x -1$. Note
that $Y$ is      a
   special $\P$-semimartingale.
   \end{definition}

\begin{remarque}
 In the case where the filtration $ \mathbb{G}$ is generated by a multidimensional Brownian motion $W$ and the jump process $N$,   the BSDEJ involves a quadratic term:
 \begin{equation*}
\left\lbrace
\begin{split}
 &dY_t=\Big[ \sum_{i=1}^d  g( \yi_t)  \lambda^i_t- U_t+\delta_t
Y_t  + \frac{1}{2} |Z_t|^2  \big] \, dt
+ Z_t .  dW_t +  \sum_{i=1}^d { \yi _t}dN^i_t\\
&Y_T= \bar{U}_T
\end{split}
\right.
\end{equation*}
Such BSDE have been studied recently in the case where the terminal condition  is bounded  and typically appear in problems from pricing-hedging derivative options by indifference pricing or/and maximization of expected exponential utility including jumps on the wealth portfolio; see for instance Becherer \cite{B06}, Morlais \cite{M09a, M09b} , Lim and Quenez \cite{LQ11}, Ankirchner, Blanchet-Scalliet and Eyraud-Loisel  \cite{ABL09}  and  Schroder and Skiadas \cite{SS08} for some recent references. However, all existence and comparison results for such equations assume that the terminal value $Y_T$ is bounded; here, we relax this condition.
In  a recent work,   El Karoui, Matoussi and Ngoupeyou \cite{ElkMN11} have  obtained an existence result of a general  class of QBSDEJs and unbounded final condition.

\end{remarque}

 We first establish a  recursion  relation  for
 solutions of \eqref{eqback} which implies the  uniqueness of the solution:
\begin{proposition} \label{Recursion} Let $(Y,M^{Y,c},\y)
\in D^{\exp}_0\times {\cal M}_{0,loc}^c \times
 \mathcal{L}^2(\lambda)
 $  be a  solution of the BSDEJ
(\ref{eqback}).  Then, $Y$ satisfies the following recursion
equality:  for any  stopping time $ \tau$ valued in $[t,T]$,
\begin{equation}\label{recursivite}
Y_{t}=-  \ln \, \Es^\P\left[\exp \Big(  - Y_\tau + \int_t^\tau
(\delta_s Y_s - U_s)ds\Big) \Big{\vert} {\cal G}_{t} \right].
\end{equation}
Moreover the BSDEJ \eqref{eqback} admits at most one solution which belongs to $ D^{\exp}_0\times {\cal M}_{0,loc}^c\times
 \mathcal{L}^2(\lambda)$.
  \end{proposition}
\noindent \textbf{Proof:}  (i) Assuming that $(Y,M^{Y,c},\y)$ is a
solution of \eqref{eqback}, we define $X_t: = Y_t-Y_0 -\int_0^t (\delta _sY_s- U_s)ds$
and $Z_t:=e^{-X_t }$.  It\^o's formula leads to $dZ_t = Z_{t-} \left[ -  dM^{Y,c}_t+\sum_{i=1}^d\left(e^{-{ \yi _t }}-1\right)dN^i_t \right]$.
Hence, $Z$ is a non-negative local martingale. Assuming that $Z$
is a martingale, one obtains, for $t<\tau<T$:
\begin{equation}
\label{exponentiel} e^{-Y_t }=\Es^\P\left[\exp \left(-Y_\tau +
\int_t^\tau (\delta_s Y_s - U_s)ds \right)\Big{\vert} {\cal
G}_t\right] \,.
\end{equation}
Otherwise, using a localizing sequence  $\tau_n$ such that  the stopped process $ Z^{\tau_n} $ is a martingale,    we obtain
\eqref{exponentiel} with  $\tau_n\wedge\tau$ instead of $\tau$. By the integrability
Assumption \ref{hyp2} and the assumption that $Y\in D_0^{\exp}$,
we obtain   a $\mathbb{P}$-integrable upper bound for the
right-hand side of \eqref{exponentiel} and letting $n$ go  to
infinity,  we obtain \eqref{recursivite}
for $\tau$  by dominated convergence. \\
(ii) \textit{Uniqueness of the solution of the  BSDE \eqref{eqback}:}   Assume  that $(Y,M^{Y,c}, \y)$ and $(\bar{Y},
{M}^{\bar Y,c}, \bar {\y})$ are  two solutions of (\ref{eqback})
in $ D^{\exp}_0\times  {\cal M}_{0,loc}^c \times \mathcal{
L}^2(\lambda) $. Suppose that, for some $t\in [0,T]$, the set
$A=\{Y_t>\bar{Y}_t\}\in {\cal G}_t$ satisfies  $\P(A)>0$ and
define $\tau=\inf\{s\ge t \vert \bar {Y}_s\ge Y_s\}$, so that
$\bar Y_\tau \ge Y_\tau$.  Since $Y_T=\bar {Y}_T$, one has $\tau
\leq T$,  and:
$$\int_t^\tau(\delta_sY_s- U_s)ds -Y_\tau>\int_t^\tau(\delta_s\bar Y_s-
U_s)ds -\bar Y_\tau \hbox{ on  } A, $$ then from the recursion relation \eqref{recursivite}, it follows that
$$\exp\left(- Y_t \right)=\Es^\P\left[\exp\left( \int_t^\tau \delta_sY_s- U_s)ds - Y_\tau\right)\Big{\vert }{\cal G}_t\right]> \exp\left(- \bar {Y}_t\right) \hbox{ on } A $$ which implies that $Y_t<\bar{Y}_t$ on $A$ in contradiction with
the definition of $A$; therefore  $Y$ and $\bar Y$ are
indistinguishable. It follows that $ M^{Y,c} = M^{\bar Y ,c} $ and $  \y =  \bar \y $.
\finproof
\begin{remarque}
In the case $\delta =0$, the process $Y$, part of the solution of
\eqref{eqback}, is given in a closed form as
$$ Y_t = -\ln \Es^\P\left[\exp \left(-\bar U_T - \int_t^T
U_s ds \right)\Big{\vert} {\cal G}_t\right] \,.$$ In the case
$ U \equiv 0$, we recognize  the dynamic  entropic risk measure
studied, in particular,   by Barrieu and El Karoui \cite{BElk08}.
\end{remarque}

 The   main result  of this section   gives  the BSDE description of the value process of our robust control problem. 
  It extends earlier works by Skiadas  \cite[Theorem 5, pp. 482]{S03}   and Bordigoni, Matoussi and Schweizer \cite[Theorem 12]{BMS07}  (see also Lazrak and Quenez \cite{LQ03} and Schroder and Skiadas \cite{SS99}).

\begin{theoreme}\label{existence}
Assume  \textbf{(A1)} and \textbf{(A2)}. Then, there exists a unique triple   $(Y,M^{Y,c}, \y) \in
D^{\exp}_0\times {\cal M}_{0}^p \times  \mathcal{L}^2(\lambda )  $ solution of \eqref{eqback}. Furthermore, the  optimal
measure $\Q^*$ solution of \eqref{objectif} admits the
Radon-Nikodym density $Z^{\Q^*}= {\cal {E}}(L)$ w.r.t. $\P$ where
\begin{equation}\label{densitel}dL_t=-dM^{Y,c}_t +\sum_{i=1}^d \left(e^{-  \yi _t
}-1\right)dN^i_t, \quad   {L_0=0}.\end{equation}
\end{theoreme}

 \noindent
\textbf{Proof:} We divide the proof in three steps. We first prove that the value process $V$ of our control problem is a $\P$-special semimartingale, i.e., $ V = V_0 + M^V + A^V$  with $ dM^V_t =d M^{V,c} _t+   \sum_{i=1}^d { v^i _t}dN^i_t$.   Secondly, we prove  that $(V,  M^{V,c}, v )$ is a  solution of the BSDE \eqref{eqback}. Finally,   we show that    $(V,M^{V,c}, v) \in
D^{\exp}_0\times {\cal M}_{0}^p\times  \mathcal{L}^2(\lambda)  $. \\[0.3cm]
\textbf{Step 1:}
We  embed the minimization of  $\Gamma (\Q)$  in a stochastic control problem and we use mainly  the martingale optimality principle   from    El Karoui  \cite{Elk82}   (Theorem 1.15,  Theorem 1.17 and Theorem 1.21) to get our result. To that end,  we
introduce a few more notation.    We define
the minimal conditional cost
$$  J (\tau, \Q) := \Q - \mbox{ess inf}_{\Q' \in \mathcal{D} (\Q, \tau)} \Gamma (\tau, \Q')$$
with $ \Gamma (\tau, \Q') := \mathbb{E}_{\Q} \left[    \mathcal U_{0,T}^{\delta} +  \mathcal R_{0,T}^{\delta} (\Q')     \, | \,  \mathcal{G}_{\tau} \right]$ and $ \mathcal{D} (\Q, \tau) = \{ Z^{\Q'}  \,|\,  \Q' \in \mathcal{Q}_f \; and  \;  \Q' = \Q \; \mbox{on}  \; \mathcal{G}_{\tau} \}$.  So, we can write our minimization problem as
$$\mbox{} \quad { \inf_{\Q \in \mathcal{Q}_f }\Gamma (\Q) =
\mathbb{E}^\P\left[ J (0,\Q)  \right]}
$$
by using the dynamic programming equation and the fact that $ \Q = \P$ on  $ \mathcal G_0$ for every $ \Q \in \mathcal Q_f$. A measure
 $ \tilde{\Q} \in \mathcal Q_f $ is called \textit{optimal} if it minimizes $\mb{ \Q }\mapsto \Gamma (\Q)$ over $ \Q \in \mathcal Q _f$.
\noindent We know from Proposition \ref{optimal} (or Theorem 9 and Theorem 12 in \cite{BMS07}) that there exists an optimal  $ \Q^*$ which belongs to  $ {\cal Q}_f^e $, hence, w.l.o.g.,  we  restrict our attention to   minimize $ \Q \longmapsto \Gamma (\Q)$   over $ \Q \in {\cal Q}_f^e $}. For each $ \Q \in \mathcal Q_f^e$ and $ \tau \in \mathcal I $, where ${\mathcal I} $ is the set of   $\mathbb G $-stopping times
valued in $[0,T]$, we define

$$  V (\tau, \Q) := \Q - \mbox{ess inf}_{\Q' \in \mathcal{D}
(\Q, \tau)}   \mathbb{E}_{\Q'} \left[    \mathcal
U_{\tau,T}^{\delta} +  \mathcal R_{\tau,T}^{\delta}     (\Q') \, |
\, \mathcal{G}_{\tau} \right]
$$
 which is the \textit{value} of the control problem started at time $\tau$  and assuming one has used the model $\Q$ up to time $\tau$. By using the Bayes formula and the definition of $ \mathcal R_{\tau,T}^{\delta}     (\Q') $, one can easily prove that $ V (\tau, \Q) =  V (\tau) $ does not depend on $\Q \in \mathcal Q_f^e$.  Moreover,  comparing the definitions of $ V(\tau)$ and $ J (\tau, \Q)$  yields   for  $\Q \in \mathcal Q_f^e$

 $$
 J (\tau, \Q) =\sd_{\tau} V (\tau) + \int_0^\tau  \sd _s U_s ds
+\int_0^{\tau}  \delta_s \sd _s\ln Z^{\Q}_sds + \sd _{\tau} \ln Z^{ \Q}_{\tau}
$$
because we can also take the ess inf for $ J(\tau, \Q)$ under $ \P \sim \Q$.
  From the martingale  optimality principle proved in    \cite{BMS07} (Proposition 13 pp.140),
there exists an adapted RCLL process $J^\Q =
(J^\Q_t)_{0\leq t \leq T}$ which is a right closed $\Q$-submartingale
such that $ J^\Q_{\tau} = J (\tau,\Q)$. Thus we can choose an adapted RCLL process $ V = (V_t)_{ 0 \leq t \leq T}$ such that
 $ V_{\tau} = V (\tau) =  V (\tau, \Q), \quad  \P-a.s.$ for $ \tau \in \mathcal I$ and
 $ \Q \in \mathcal Q_f ^e$, and then we get,  for each $ \Q \in \mathcal Q_f^e$,
  \begin{equation}
  \label{submartingale}
  J^{\Q}=\sd V+ \int_0 \sd _s U_s ds+\int_0 \delta_s \sd _s\ln Z^{\Q}_sds + \sd  \ln Z^{ \Q}.
  \end{equation}
As $ \P \in \mathcal Q_f^e $ and $ J^\P$ is a $\P$-submartingale   (from Proposition 13  pp. 140 in \cite{BMS07}), equation
\eqref{submartingale}  yields that $J^\P =  \sd V+ \int \sd _s U_s ds$. Thus $ V$  is a $\P$-special semimartingale, i.e., its canonical decomposition can be written as
$$
V = V_0 + M^V + A^V.
$$
Since $ \sd$ is uniformly bounded from below and $ J^\P$ is a $ \P$-submartingale,   Assumption  \textbf{(A2)} implies that $ M^V$ is a true $\P$-martingale and that
$ dM^V_t = dM^{V,c}_t +   \sum_{i=1}^d { v^i_t  }dN^i_t $ where $ M^{V,c}$ is  a continuous $\P$-martingale.\\\\
\noindent
\textbf{Step 2:}   We  now prove  that $(V,M^{V,c}, v )$ is  solution of the BSDEJ \eqref{eqback} where $ v := (v^1, \cdots, v^d)$.\\\\
For  ${\Q}
\in {\cal Q}^{e}_f$,  we denote by $L^{\Q}$  the stochastic logarithm of $Z^\Q$, i.e., the $\mathbb P$-local
martingale such that $ dZ^{\Q}_t=Z^{ \Q} _{t-} dL^{\Q}_t$.

From
Assumption  \ref{hyp},  the $\P$-local martingale  $L^{\Q}$ admits the
decomposition $dL^{\Q}_t =dL^{{\Q},c}_t +\sum_{i=1}^d  \lqi_t  dN^i_t $,
where $L^{\Q,c}$ is a continuous $\P$-local  martingale, and
$\lqi$ are predictable processes, and one has
\begin{equation}\label{Lndynamics} d\ln Z^{\Q}_t=dL^{{\Q},c}_t-{1\over
2}d{\langle L^{\Q,c}\rangle}_t +\sum_{i=1}^d  \ln(1+\lqi_t)
dN^i_t+\sum_{i=1}^d (\ln(1+\lqi_t) -\lqi_t)\lambda^i_tdt.
\end{equation}
\noindent
Using integration by parts formula,  we obtain after some simple computations and using \eqref{submartingale} and
\eqref{Lndynamics}:

\begin{eqnarray*}
dJ^{\Q}_t &=&\sd _t\left((-\delta_tV_t+ U_t)dt+(dV_t+ d\ln
Z^{\Q}_t)\right) \label{JLdynamics} \\ &=&\sd
_t\Big[(-\delta_tV_t+ U_t)dt+dM^{V,c}_t+dA^V_t+  d L^{\Q,c}
_t-{1\over 2}d{\langle L^{\Q,c}\rangle}_t \nonumber
\\&&+\sum_{i=1}^d (\vi_t + \ln(1+\lqi_t)) dN^i_t+\sum_{i=1}^d
(\ln(1+\lqi_t) -\lqi_t)\lambda^i_tdt  \Big]\, \nonumber
\end{eqnarray*}
\noindent From Girsanov's theorem, the   processes $ { \widetilde N^i  } $ and
${ \widetilde M^{ c}  } $ defined as $d\widetilde N^i_t=dN^i_t- \lqi_t \lambda^i_tdt, \hbox{ and }
d\widetilde M^{c}_t=d(M^{V,c}_t+  L^{\Q,c}_t)-d\langle M^{V,c}+
L^{\Q,c},L^{\Q,c}\rangle_t$  are ${\Q }$-local martingales, and:
 \begin{eqnarray*}dJ^{\Q}_t&=&\sd _t\Big[(-\delta_tV_t+ U_t)dt+d\widetilde M^{c}_t+dA^V_t+  d\langle M^{V,c}+
L^{\Q,c},L^{\Q,c}\rangle_t -{1\over 2}d{\langle L^{\Q,c}\rangle}_t
\\&&+\sum_{i=1}^d (\vi_t + \ln(1+\lqi_t)) d\widetilde
N^i_t+\sum_{i=1}^d \left(  \lqi _t(\vi _t -1) + (1+\lqi _t)
\ln(1+\lqi_t)\right)
 \lambda^i_tdt  \Big].
\end{eqnarray*}
In order that the process $J^{\Q}$ is a ${\Q}$-submartingale for
each ${\Q} \in {\cal Q}^{e}_f$, we impose that   its finite
variation part is a non-decreasing process.
\begin{eqnarray}A^V_t&=& - {\rm ess}\, \inf_{{\cal Q}^e_f}\int_0^t
( U_s-\delta_sV_s)ds + \langle M^{V,c}+ L
^{\Q,c},L^{\Q,c}\rangle_t -{1\over 2} {\langle L^{\Q,c}\rangle}_t
   \nonumber \\&&+\sum_{i=1}^d\int_0^t  \left(  \lqi _s (\vi _s -1) + (1+\lqi _s)
\ln(1+\lqi_s)\right)
 \lambda^i_sds. \label{A}\end{eqnarray}
 To find the $ess\inf$, we divide \eqref{A} in two parts, the continuous part and the discontinuous
 part;
hence we have two optimization problems:
\begin{eqnarray*}
 \label{A_suite}
A^V_t&=&\int_0^t(\delta_sV_s- U_s)ds- {\rm ess}\inf_{{\cal
Q}^e_f}\{\langle M^{V,c},L^{\Q,c}\rangle_t +{1\over 2}\langle
L^{\Q,c}\rangle_t\}  \\& - &{\rm ess}\inf_{{\cal
Q}^e_f}\sum_{i=1}^d \int_0^t  \left(  \lqi_s (\vi_s -1) +
(1+\lqi_s) \ln(1+\lqi_s)\right)
 \lambda^i_sds.
\end{eqnarray*}
 It is proved  in \cite{BMS07} that the first  infimum is obtained for   $L^{{\Q},c}= -
M^{V,c}$ and $-{\rm ess}\inf_{{\cal Q}^e_f}\{\langle
M^{V,c},L^{\Q,c}\rangle+{1\over 2}\langle L^{\Q,c}\rangle
\}={1\over 2 }\langle M^{V,c}\rangle$. The second part of the optimisation problem reduces to find
the optimal $\lqi$, solution of  ${\rm ess}\inf \left(  \lqi_s (\vi_s -1) + (1+\lqi_s)
\ln(1+\lqi_s)\right)$ which is an easy task, the solution being  $\lqsi_s= e^{- \vi
_s}-1$, which leads to $$-{\rm ess}\inf \left(  \lqi_s (\vi_s -1) + (1+\lqi_s)\ln(1+\lqi_s)\right)\\=e^{- \vi_s} +\vi_s-1  =  g( \vi_s)$$ where $g (x)=   e^{ -x}+x-1  $. Therefore,
\begin{equation*}\label{A_fin1}
A^V_t=\int_0^t(\delta_sV_s- U_s)ds+ {1\over 2 }\langle
M^{V,c}\rangle_t + \int_0^t \sum_{i=1}^d g(\vi_s)\lambda^i_sds\,.
\end{equation*}
\noindent It follows that $(V,M^{V,c},  v )$ is a solution of
\begin{equation}\label{equa}
dV_t=\left(\delta_t V_t- U_t+\sum_{i=1}^d g (\vi
 _t)\lambda^i_t\right)dt +\frac 12 d\langle M^{V,c}\rangle_t + dM^{V,c}_t+\sum_{i=1}^d \vi_tdN^i_t,\,
V_T= \bar{U}_T
\end{equation}
 \noindent  Furthermore  there exists a  solution of the  BSDEJ  \eqref{eqback}  and  the optimal probability measure $\Q^*$ is characterized by
 its Radon-Nikodym density
$$dZ^{\Q^*}_t = Z^{\Q^*}_{t^-}dL _t,\qquad  dL _t=
- dM^{V,c}_t +\sum_{i=1}^d\left(e^{-\vi_t }-1\right)dN^i_t.$$
\noindent
\textbf{Step 3:}  In this step we prove that the solution
$(Y,M^{Y,c},\y)$ of the BSDEJ  \eqref{eqback} belongs to the
required spaces.   \\
From Lemma 19 and Proposition 20 in  \cite{BMS07}, we have that $Y$
 belongs to $D^{exp}_0$. To study  the space of the process  $M^{Y,c}$, we introduce the
$\P$-martingale:
$$K_t:=\Es^{\P}\left[\exp \left( \int_0^T(\delta_s Y_s
- U_s)ds -  \bar U_T\right)\Big{\vert}\cG_t\right]\,.$$ Using the
fact that $Y \in D^{exp}_0$, we obtain that the process  $K $
belongs to ${\cal M}^p $. Now, the recursive property
leads to 
$  K_t =\exp\left(- Y_t+ \int_0^t (\delta_sY_s- U_s)ds\right)\,$
and from It\^o's formula and the
canonical decomposition of $Y$,  \begin{equation}
\label{mk}dM^{Y,c}_t=-{d K^c_t\over K_{t-}}\,. \end{equation}
From
Assumption \ref{hyp}, there exists   $k^i$ and $M^{K,c}$ such that
  $K_t=K_0 +M^{K,c}_t+\sum_{i=1}^d \int_0^t    k^i_s dN^i_s$.
Hence, from \eqref{mk}:
\begin{equation*}
\begin{split}
\langle M^{Y,c}\rangle_T  & \le \int_0^T {1\over K^2_t}d\langle
K^c\rangle_t   \le \langle K^c\rangle_T \sup_{0\le t\le T}{1\over
K^2_t}\\
& \le \langle K^c\rangle_T \exp\left({2}\sup_{0\le t\le T}|Y_t|(1+||\delta||_{\infty}T)+{2 }\int_0^T|U_s|ds\right)\,.\\
\end{split}
\end{equation*}
By BDG's inequalities, there exists a constant $C$ such that for
every $p\in [1,+\infty)$:
\begin{equation}\label{BDG}
\Es^{\P}{\left[\langle K^{c}\rangle_T+\int_0^T{(k^i_t)}^2 dH^i_t\right]}^{p\over 2}\le C \Es^{\P}(\sup_{0\le t\le T}{|K_t|}^p)
\end{equation}
\noindent Since $K\in {\cal M}^p $,  we   conclude that $M^{Y,c}$
lies in the   space ${\cal M}_0^p $ for every $p \in
[1,+\infty)$. We conclude, using again BDG's inequalities. Finally, let now characterize the space of   process $y$. Using the recursive relation and the
decomposition of the process $K$ we get $\ln(K_{t^-}+k^i_t)-\ln(K_{t^-})=-y^i_t$, hence:
 $$\Es^\P\left[\int_0^T {(e^{-y^i_t}-1)}^2 dH^i_t\right]^{p\over 2}=\Es^\P\left[\int_0^T {\left({k^i_t\over K_{t^-}}\right)}^2 dH^i_t\right]^{p\over 2}\le \Es^\P\left[\left(\sup_{0\le t\le T}{1\over K^p_t}\right){\left(\int_0^T {(k^i_t)}^2 dH^i_t\right)}^{p\over 2}\right]\,.$$
\noindent Since $\displaystyle \sup_{0\le t\le T} \Big({1\over K_t}\Big)
\in \mathbf{L}^p(\P)$ for any $p\in [1,+\infty]$,
 using ($\ref{BDG}$) and Cauchy inequalities, we conclude:
\begin{equation}
\label{12}\Es^\P\left[\int_0^T {(e^{-y^i_t}-1)}^2
dH^i_t\right]^{p\over 2}<\infty.\end{equation}
 In particular
\begin{equation}
\label{estimate1:y} \Es^\P\left[\int_0^T
{(e^{-y^i_t}-1)}^2\lambda^i_tdt\right] =\Es^\P\left[\int_0^T
{\left({k^i_t\over K_{t^-}}\right)}^2 dH^i_t\right] \le \Es
^\P\left[\left(\sup_{0\le t\le T}{1\over
K^2_t}\right)\int_0^T{(k^i_t)}^2dH^i_t\right]<\infty\,.
\end{equation}
\noindent By using similar arguments, one proves that:
\begin{equation}
\label{estimate2:y}\Es^\P\left[\int_0^T {(e^{y^i_t}-1)}^2\lambda^i_tdt\right] <\infty.
\end{equation}
\noindent Moreover, by using the  inequality  $ |y |^2 \leq 2 \big(|e^{-y} -1|^2 + |e^{y} -1|^2\big), \; \; \forall  y\in \mathbb{R}$ and  \eqref{estimate1:y}-\eqref{estimate2:y} we conclude that the process $ y$  belongs to $\mathcal{L}^2 (\lambda )$.
\finproof
\begin{remarque} The martingale part of the  BSDE
solution, i.e.,  $M=-M^{Y,c}+\sum_{i=1}^d \int_0^.
(e^{-y^i_t}-1)dN^i_t $ belongs to ${\cal M}_0^p $ for any
$p\in[1,+\infty)$. Indeed, since $M^{Y,c}\in{\cal M}_0^p $, and  (\ref{12})
$$\Es ^\P{\left[\langle M^{Y,c}\rangle_T+\sum_{i=1}^d \int_0^T  {(e^{-y^i_t}-1)}^2 dH^i_t\right]}^{p\over 2}<\infty,$$
and using BDG inequality, we obtain $\Es^\P\left(\sup_{0\le t\le T}|M_t|^p \right)<\infty.$
\end{remarque}

%

\subsection{Properties of the value process}
\noindent In this part, we establish a comparison theorem  for the   class of  BSDEJs \eqref{eqback}  which is a key point  to derive the dynamic maximum principle  for the maximization problem. 

\begin{definition} \label{ordre}For two random variables $X$
and $Y$, we   write $X\leq Y$ for $X\leq Y$ a.s. For two
processes $A$ and $B$, we   write $A\leq B$ for $A_t\leq
B_t,\forall t\in [0,T], a.s.$.  We   write  $(X,A)\leq(Y,B)$
if $X\leq Y$   and $A\leq B$.
\end{definition}
\begin{theoreme}\label{comparaisonth} Assume that for $k=1,2$, $(Y^k,M^{ k,c},\yk)$ is the solution of the
BSDE (\ref{eqback}) associated with ($U^k,\bar U^k_T$). We denote
$\yud :=Y^1-Y^2$, $ \uudt:=U^1-U^2$ and  $ \uudb_T:=\bar
U^1_T-\bar U^2_T$. Then,
\begin{equation}\label{comparison}
 S_t^{\delta} \,  \yud_t   \, \le \, \mathbb{E}^{\mathbb{Q}^{*,2}}\left[\int_t^T
S_s^{\delta}  \uudt_s   ds + S_T^{\delta}   \uudb_T \Big{\vert}
{\cal G}_t\right]\,,
\end{equation}
 where  $\Q^{*,2}$  is the solution of ${\cal {P}}(U^2,
\bar U_T^2)$, i.e., the probability measure  equivalent to $\P$
with Radon Nikodym density $  Z^{\Q^{*,2}} $  given by
\begin{equation}\label{z2} dZ^{\Q^{*,2}}_t= Z^{\Q^{*,2}}_{t^-} \left(-dM^{ 2,c}_t
+\sum_{i=1}^d \Big(e^{-
\yideux_t}-1\Big)dN^i_t\right).\end{equation} In particular, if
$(U^1,\bar U^1_T)\le (U^2,\bar U^2_T)$, one obtains $Y^1_t\le Y^2_t, \, dP\otimes dt\textbf{-}a.e$.
\end{theoreme}
 \noindent
\textbf{Proof}: We denote $\y^{i,12} :=\yiun -\yideux $ and $M^{12,c}:=M^{1,c}-M^{2,c}$, then we find that:
\begin{equation*}
\begin{split}
\yud_t=& \;  \uudb_T  +\int_t^T \left(    \uudt _s
-\delta_s\yud_s\right)ds-\sum_{i=1}^d\int_t^T  \y^{i,12}_s dN^i_s
- \sum_{i=1}^d\int_t^T  \left[g(\yiun_s) -g(\yideux_s )
\right]\lambda^i_sds\\&
 +{1\over 2  } \int_t^T \left( d{\langle M^{ 2,c}\rangle}_s-  d{\langle M^{
 1,c}\rangle}_s\right)
- \int_t^T d M^{ 12,c}_s\,.
\end{split}
\end{equation*}
\noindent Note that, since  $M^{k,c}$ are continuous martingales,
\begin{equation}\label{util}
 - \langle M^{2,c}, M^{12,c}\rangle  - \frac 12 \langle M^{2,c}\rangle
  +
  \frac 12  \langle M^{1,c}\rangle  = \frac 12\langle M^{12,c} \rangle \,.
\end{equation}
\noindent Using the fact that the process $\langle M^{12,c} \rangle$ is
increasing and that   the function $g $  is convex, we get:
\begin{equation*}
\begin{split}
\yud_t&\le    \uudb_T  + \int_t^T \left(
\uudt_s-\delta_s\yud_s\right)ds +\sum_{i=1}^d
 \int_t^T   (e^{-\yideux_s }-1) \y^{i,12}_s \lambda^i_sds
  + \int_t^T d{\langle M^{ 2,c} ,M^{ 12,c}\rangle }_s\\&  -\int_t^T d M^{ 12,c}_s  -\sum_{i=1}^d
  \int_t^T
  \y^{i,12}_sdN^i_s.\\
\end{split}
\end{equation*}
 Let $N^*$ and $M^{*,c}$ be  the $\Q^{*,2}$-martingales
obtained by Girsanov's transformation from  $N$ and $M^{12,c}$,
where $d\Q^{*,2}= Z^{\Q^{*,2}}d\P$ and  where  $ Z^{\Q^{*,2} }$ is
given by
 (\ref{z2}). Then:  \begin{equation*}
\yud_t\le    \,\uudb_T +\int_t^T \Big(
 \uudt_s
 -\delta_s\yud_s\Big)ds-\sum_{i=1}^d\int_t^T \y^{i,12}_s dN^{i*}_s
-\int_t^T dM^{*,c}_s
\end{equation*}
which implies that $\yud_t\le
\mathbb{E}^{\Q^{*,2}}\Big[\int_t^T  e^{-\int_t^s \delta_r dr}
 \uudt_s ds +  e^{-\int_t^T \delta_r dr}  \uudb_T  \Big{\vert} {\cal G}_t\Big]
$. In particular,   if $(U^1, \bar U^1_T)\le( U^2 ,\bar U^2_T)$, then
$Y^1_t\le Y^2_t \hspace{0.2cm} d\P\otimes dt\textbf{-}a.e.$
$\hfill \Box$ \\\\
We  now prove standard a priori estimates for the solution of BSDEJ \eqref{eqback}.
\begin{proposition}\label{estimation}(A priori estimates)
 Let $(Y^k,M^{k,c},\yk)$  be  the solution associated
 with
$(U^k,\bar U^k_T)$  for $k=1,2$  where we assume that $(U^1, \bar
U^1_T)\le (U^2,  \bar U^2_T)$; then, there exists a constant $C>0$
such that:
\begin{equation}\label{priori}
\mathbb{E}^{\Q^{*,2}}\Big[\sup_{0\le t\le T}{|Y^{12}_t
|}^2+{\langle M^{ 12,c} \rangle }_T+ \sum_{i=1}^d\int_0^T
{|\y_t^{i,12 } |}^2 \lambda^{i,*}_t dt\Big] \le C \,
\mathbb{E}^{\Q^{*,2}}\Big[ {|\bar U^{12}_T |}^2+\int_0^T
{|U^{12}_t |}^2 dt\Big]
 \end{equation}
where  $  \lambda ^{i,*}  = \lambda ^i e^{-y^{i,2}} $ is the intensity process of $H^i$
under the probability $\Q^{*, 2}$. In the case $(U^2 , \bar U_T^2)\leq  (U^1, \bar U^1_T)$,  the same
inequality holds with $ \Q^{*,1}$.
\end{proposition}
\noindent \textbf{Proof: } \label{Appendix:1}   Using It\^o's
formula:
\begin{eqnarray*}
d({\yud_t})^2&=&2\yud_t\left[\left(\delta_t\yud_t-\uudt_t\right)dt+
\frac 12 d\langle M^{ 1,c}\rangle_t- \frac 12d\langle M^{ 2,c}
\rangle_t\right]+d\langle M^{ 12,c} \rangle_t \\&
+&2\yud_t\left[\sum_{i=1}^d  \left(g( \yiun_t) - g(\yideux_t)
\right)\lambda^i_t\right]dt+\sum_{i=1}^d ( \y^{i,12}_t
)^2\lambda^i_t dt + d\mbox{mart} _t
\end{eqnarray*}
where   $d\mbox{mart} _t= 2\yud_{t-} \left[d M^{ 12,c}_t
+\sum_{i=1}^d \y^{i,12}_t dN^i_t\right]+ \sum_{i=1}^d
\left(\y^{i,12}_t \right) ^2dN^i_t$   corresponds to  a local
martingale. Assuming
$(U^1, \bar U^1_T)\leq (U^2, \bar U^2_T)$, it follows from the
comparison Theorem  \ref{comparaisonth} that $Y^1\le Y^2$. Using
the   relation (\ref{util})  and the convexity property of the
function $g $, we get

%
 \begin{eqnarray*}
{(\yud_t)}^2&+&\int_t^T d\langle M^{ 12,c} \rangle_s
 \le   {(\uudb_T)}^2+2\int_t^T \yud_s
\left[-\delta_s \yud_s+ \uudt_s \right]ds\\&+&2\sum_{i=1}^d  \int_t^T
\yud_s\,(e^{{-\yideux_s}}-1)\y^{i,12}_s \lambda^i_s ds
 -2  \int_t^T \yud_sd\langle M^{ 1,c},M^{ 2,c}\rangle_s\\&+
 &\int_t^T \yud_sd\langle M^{ 2,c}\rangle_s-\sum_{i=1}^d \int_t^T  \left(\y^{i,12}_s\right)  ^2 \lambda
^i_s ds+ \int_t^T d\mbox{mart}_s \,.
\end{eqnarray*}
Hence, we  obtain the following inequality:
\begin{eqnarray}
 {(\yud_t)}^2+\int_t^T d\langle M^{ 12,c} \rangle_s
  &\le & {(\uudb_T)}^2+2\int_t^T \yud_s\left[-\delta_s\yud_s+ \uudt_s \right]ds -\sum_{i=1}^d
  \int_t^T
{\left(\y^{i,12}_s\right)}^2 \lambda ^{i*}_sds \nonumber\\& +&\int_t^T
d\mbox{mart}^*_s\,, \label{inegalite22}
\end{eqnarray}
\noindent where $\mbox{mart}^* $ is a $\Q^{*,2}$ local martingale and $\lambda
^{i*}_s:= \lambda _s^i e^{-\yideux_s}$ is the intensity of $H^i$
under $\Q^{*,2}$. From the obvious inequality
\begin{equation*}
(\yud_t)^2-{1\over \delta_t}(\uudt_t)\yud_t\ge {-1\over 4
\delta^2_t}{(\uudt_t) }^2
\end{equation*}
\noindent and the positivity of $\delta$, we deduce easily
that
\begin{equation}\label{use}
-\yud_t  \left(\delta_t  \yud_t  -   \uudt_t \right)\le { 1\over 4
\delta_t}{(\uudt_t)}^2
\end{equation}
Plotting relation (\ref{use}) in (\ref{inegalite22}) and   using the fact that    the
process $\delta$  is bounded below,   there exists a constant
$C>0$ such that:
\begin{eqnarray*}\label{finalesti}
 \mathbb{E}^{\mathbb{Q}^{*,2}}\left[\sup_{t\in
[0,T]}{|\yud_t|}^2+\langle M^{ 12,c} \rangle_T+\sum_{i=1}^d\int_0^T
 {|\y^{i,12}_t |}^2 \lambda^{i,*}_t dt\right]
 \le C
\mathbb{E}^{\mathbb{Q}^{*,2}}\left[{|\uudb_T|}^2+\int_0^T
{|\uudt_t|}^2dt\right].
 \end{eqnarray*}
 Permuting $Y^1$ and
$Y^2$  and assuming $(U^1,\bar U^1_T)\geq (U^2 , \bar U^2_T)$
leads to similar inequality. $\hfill\Box$ \\\\
\noindent
As a direct  consequence of the comparison theorem, we prove the concavity property for the BSDE solution.
\begin{theoreme}\label{concave}(Concavity property)
Define the map $F:D^{\exp}_1\times L^{\exp}\longrightarrow
D^{\exp}_0$ as $$F(U,\bar U)=V$$ where $(V,M^{V,c},\vh)$ is the
solution associated with $(U,\bar U)$. Then $F$ is concave,
namely, for all $\theta \in (0,1)$ and $(U^1,\bar U^1_T),(U^2,\bar
U^2_T) \in D^{\exp}_1\times L^{\exp}$ :
 $$F\left(\theta U^1+(1-\theta)U^2,\theta \bar U^1_T+(1-\theta)\bar U^2_T\right)\ge \theta F(U^1,\bar U^1_T)+(1-\theta)F(U^2,\bar U^2_T)
.$$
\end{theoreme}
\vspace{2mm}
 \noindent
\textbf{Proof:}
 Let $(V^k,M^{k,c},\vh^k )$  be  the solution of BSDE \eqref{eqback}  associated with $(U^k,\bar U^k_T)
 \in D^{\exp}_1\times L^{\exp} $, then for any $\theta \in (0,1)$:
 \begin{equation}\label{Vpre}
 \begin{split}
 d(\theta V^1_t-(1-\theta)V^2_t) &=\left[\delta_t(\theta
V^1_t+(1-\theta)V^2_t)-(\theta U^1_t+(1-\theta)U^2_t)\right]dt
 \\&+\,\theta {d{\langle M^{ 1,c
} \rangle}_t}+(1-\theta){d{\langle M^{ 2,c} \rangle}_t}
+d\big(\theta M^{ 1,c}_t+(1-\theta)M^{ 2,c}_t\big)
\nonumber\\&+\,\sum_{i=1}^d \big[\theta \vh_t^{1,i} +(1-\theta
)\vh_t^{2,i}\big]\, dN^i_t  + \sum_{i=1}^d\big[\theta\,
g(\vh_t^{1,i})+
 (1-\theta)g(\vh_t^{2,i})\big]
\lambda^i_t dt.
  \end{split}
  \end{equation}
 We recall the following general result: Let $X$ and $Y$ be two continuous
martingales. Then, for all $\theta \in (0,1)$, $ \theta \langle
X\rangle+(1-\theta)\langle Y\rangle-\langle\theta
X+(1-\theta)Y\rangle $ is an increasing process. Indeed, we have:
 \begin{eqnarray*}
&&\langle \theta X+(1-\theta)Y\rangle-\theta \langle X\rangle
-(1-\theta)\langle Y\rangle
\\&=&(\theta^2-\theta)\langle
X\rangle+\big({(1-\theta)}^2-(1-\theta)\big)\langle
Y\rangle+2\theta(1-\theta)\langle X,Y\rangle\\&
=&\theta(\theta-1)\big[\langle X\rangle+\langle Y\rangle-2\langle
X,Y\rangle\big]  =\theta(\theta-1) \langle X-Y\rangle.
\end{eqnarray*}
 Therefore,  using the convexity property of the function
$g$, we get:
\begin{eqnarray}\label{CONCA}
\theta V^1_t+(1-\theta)V^2_t &\le & \big(\theta \bar
U^1_T+(1-\theta) \bar U^2_T\big)-\int_t^T\big[\delta_s(\theta
V^1_s+(1-\theta)V^2_s)-(\theta U^1_s+(1-\theta)U^2_s)\big]ds
\nonumber\\ &&-\int_t^T {d{\langle \theta M^{ 1,c}+(1-\theta)M^{
2,c}\rangle}_s}
-\int_t^T d(\theta M^{ 1,c}_s+(1-\theta)M^{ 2,c}_s)\nonumber\\
\lefteqn{-\sum_{i=1}^d \int_t^T (\theta \vh_s^{1,i}+(1-\theta
)\mb{\vh_s^{2,i}})dN^i_s-\sum_{i=1}^d \int_t^T  g({\theta
\vh_s^{1,i}+(1-\theta)\vh_s^{2,i}})
 \lambda^i_s ds.}
\end{eqnarray}
Let  $(V^\theta,M^{ \theta,c},\vh^{\theta})$  be  the solution of
the BSDE associated with $\big(\theta U^1+(1-\theta)U^2,\theta
\bar U^1$ $+(1-\theta)\bar U^2\big)$ and set  $M^{V,c,
\theta}=\theta M^{ 1,c}+(1-\theta)M^{ 2,c}$ and for $i=1,\dots,
d$, $\widehat{v}^{\theta,i} =\theta
\vh^{1,i}+(1-\theta)\vh^{2,i}$.  Then,  using (\ref{CONCA}):
\begin{eqnarray*}
 \theta V^1_t&+&(1-\theta )V^2_t -V_t^{\theta} \le  \int_t^T
\delta_s(V_s^{\theta}-(\theta
V^1_s+(1-\theta)V^2_s)) \, ds -\int_t^T {d{\langle M^{V,c, \theta}\rangle}_s}+\int_t^T {d{\langle M^{ \theta,c}\rangle}_s} \nonumber\\
 &&-\int_t^T
d(M_s^{V,c, \theta}-M^{ \theta,c}_s) +\sum_{i=1}^d \int_t^T
\big(g( {\vh_s^{\theta,i}} )-g ({\widehat{v}_s^{\theta,i}}
)\big)\lambda^i_s ds -\sum_{i=1}^d\int_t^T  (
\widehat{v}_s^{\theta,i}-\vh_s^{\theta,i})dN^i_s.
\end{eqnarray*}
Using $\eqref{util}$ and the  convexity  property of the function
$g$, we get:
\begin{eqnarray*}\label{Concademo1}
\theta V^1_t&+&(1-\theta )V^2_t -V_t^{\theta} \le
\int_t^T\left[\delta_s(V_s^{\theta}-(\theta
V^1_s+(1-\theta)V^2_s)\right]\, ds\\& +&\sum_{i=1}^d\int_t^T \left(
e^{-{\vh_s^{\theta,i}}}-1\right)(\widehat{v}_s^{\theta,i}
-\vh_s^{\theta,i})\lambda^i_s ds - \int_t^T d({\langle M^{
\theta,c},M^{V,c, \theta} \rangle}_s +{\langle M^{
\theta,c}\rangle}_s)  \\&&-\int_t^T d(M_s^{V,c, \theta}-M^{
\theta,c}_s) -\sum_{i=1}^d\int_t^T   ( \widehat{v}_s^{\theta,i} -
\vh^{\theta,i}_s)dN^i_s.
\end{eqnarray*}
\noindent Therefore, we find the following inequality:
\begin{eqnarray*}
\theta V^1_t&+&(1-\theta )V^2_t -V_t^{\theta}  \leq \int_t^T\big[\delta_s(V^\theta_s-(\theta
V^1_s+(1-\theta)V^2_s)\big]\, ds\\&-&\sum_{i=1}^d \int_t^T
 \big(
 \widehat{v}_s^{\theta,i}
-\vh_s^{\theta,i})(dN^i_s-(e^{-{\vh_s^{\theta,i}}}-1)\lambda^i_sds\big)
\\& &-\int_t^T d \left((M_s^{V,c, \theta}-M_s^{ \theta,c})+{\langle
M^{V,c, \theta}-M^{ \theta,c},{M^{ \theta,c}}\rangle_s}\right).
 \end{eqnarray*}
\noindent Let $\Q^{*,\theta}$ be the probability measure equivalent to $\P$
with  Radon-Nikodym  density  given by $$ dZ^{\Q^{*,\theta}}_t= Z^{\Q^{*,\theta}}_{t^-}
\left[-{dM^{ \theta,c}_t}+\sum_{i=1}^d
(e^{-\vh_t^{\theta,i}}-1)dN^i_t\right], $$ then using integration
by parts and  Girsanov's theorem, taking $\Q^{*,\theta}$-conditional expectations, we have
$ \sd _t\left(\theta V^1_t+(1-\theta )V^2_t  -V^\theta_t\right) \le 0$, which gives the result.
$\hfill\Box$

\section{The second optimization problem}
 \noindent In this section, we assume that $U_s= U(c_s)$ and $\bar U_T=\bar U(\psi)$ where $U$ and
 $\bar U$ are given utility functions, $c$ is a given  non-negative $\gg$-adapted
 process and $\psi$ a $\cG_T$-measurable non-negative random
 variable. We fix  a probability $\widetilde \P$
equivalent to $\P$ with a Radon-Nikodym density $\wt Z$  with
respect to $\P$ of the from :
\begin{equation}\label{risk}
d\wt Z_t =  \wt Z_{t-} (\theta_tdM^c_t+ \sum_{i=1}^d
(e^{-z^i_t}-1)dN^i_t),\, \wt Z_0=1\;.
 \end{equation}
\subsection{Formulation of the problem}
 \noindent We study  the following  optimization problem  of  the robust maximization initial problem \eqref{pb}:
\begin{equation*}
\begin{split}
&\sup_{(c,\psi)\in {\cal A}(x)}\Es^{\Q^*}\left[ \int_0^T
\sd _s{U(c_s)}ds+ \sd _T \bar U(\psi)
\right]+\Es^{\Q^*}\left[\int_0^T \delta_s \sd _s \ln Z^{\Q^*}_s
ds+\sd _T\ln Z^{\Q^*}_T\right]\\& =:\sup_{(c,\psi)\in {\cal A}(x)}
V^{x,\psi,c}_0
\end{split}
\end{equation*}
\noindent where ${\cal A}(x)$  is the  set \mb{of} admissible control parameters,   $V_0$ is the   value at initial time of the value process
$V$, part of the solution $(V,M^{V,c}, v)$  of the BSDEJ
(\ref{eqback}) in the case $U_s=U(c_s)$ and $\bar U_T=\bar
U(\psi)$. Here, $\Q^*$ is the optimal
 measure for ${\cal {P}}( U(c), \bar U(\psi))$, and depends
on $(c,\psi)$. The preferences are modeled by the utility functions $U$ and $ \bar{U}$ which
satisfy the following conditions:
\begin{Hypothese} \label{utility} The utility functions $U$ and $\bar{U}$
  satisfy the usual   conditions:\vspace{0.2cm}\newline i) Strictly
increasing and concave. \vspace{0.1cm}\newline ii) Continuous
differentiable on the set $\{U>-\infty\}$ and  $\{\bar
U>-\infty\}$, respectively, \vspace{0.1cm}\newline iii)
$U'(\infty):=\lim_{x\rightarrow \infty}U'(x)=0$ and $\bar
U'(\infty):=\lim_{x\rightarrow \infty}\bar U'(x)=0$,
\vspace{0.1cm}\newline iv) $U'(0):=\lim_{x\rightarrow
0}U'(x)=+\infty$ and $U'(0):=\lim_{x\rightarrow 0}\bar
U'(x)=+\infty$, \vspace{0.1cm}\newline
v) Asymptotic elasticity $AE(U):=\displaystyle \lim\sup_{x\rightarrow +\infty}{xU'(x)\over U(x)}<1$.\\
\end{Hypothese}
  \begin{definition}   ${\cal A}(x)$ is the set of control parameters
$(c,\psi) \in {\cal H}^2([0,T])\times
{\mathbf{L}}^2(\Omega,{\cal{G}}_T )$
such that
\begin{equation}
\label{constraint}
 \Es^{\widetilde
{\mathbb{P}}}\big[\int_0^T c_t dt+\psi \big]\le x\,,
\end{equation}
 and $(U(c
),\bar U(\psi)) \in D^{\exp}_1\times L^{\exp} $ and $(\bar c
U'(c)$,$\bar \psi \bar U'(\psi))$ $\in D^{\exp}_1\times L^{exp}$
for any pair $(\bar c,\bar \psi) \in {\cal H}^2([0,T])\times
{\mathbf{L}}^2(\Omega,{\cal{G}}_T )$, and the process $\exp{(\gamma \int_0^{\cdot} |U(c_t)|dt)}$  (resp. $ \exp{(\gamma \int_0^{\cdot} |c_t||U'(c_t)|dt)}) $ belongs to the class $\bf [D]\rm$ (see Dellacherie and Meyer, pp.89, Chapter VI \cite{DM80} for definition).
\end{definition}

\begin{remarque}
Under our assumptions,  the set $ \mathcal A (x)$  is convex  and closed in the topology of convergence in measure (see, Cuoco \cite{C97} Lemma B3., pp.70).
\end{remarque}
\noindent In order to clarify and motivate the constraint \eqref{constraint}  satisfied  by the control parameters $ (c, \psi)$, we  present a generic example in a financial market where the process  $c$ can be interpreted as  a consumption and $\psi$ as a terminal wealth :
\begin{exemple}[\textit{Consumption-investment problem}]
\label{example2}
We assume the same model as in Example \ref{exemple1}, and we consider  a  financial market consisting of $d+p +1$ assets. The
\textit{savings account}  is assumed to be constant equal to 1, the  prices of the $d+p$
\emph{risky assets} are $\gg$-semi-martingales  given by
\begin{equation}
dS^i_t=S^i_{t^-}\left[\mu^i_tdt+\sum_{j=1}^d
\varphi^{i,j}_tdN^j_t+ \sum_{k= 1}^{
p}\sigma^{i,k}_tdW^{k}_t\right], i=1,\dots, d+p
\end{equation}
where 
$\sigma$ is
a $(d+p) \times d $ volatility matrix $(\sigma^{i,k}, i=1,\dots,
d+p;k=1,\dots, p)$ and  $\varphi$   is a $(d+p) \times d$ matrix
$(\varphi^{i,j}, i=1,\dots, d+p;j=1,\dots, d)$.   We   note
$\Sigma$ the $(d+p) \times (d +p )$ matrix  $ \Sigma =
[\sigma, \lambda \varphi]$, where $\lambda \varphi$ is the matrix with coefficients ($\lambda ^i\varphi^{i,j}) $.\\
  Given an initial wealth $x$ and a policy $(c,\pi)$, the wealth process  $ (X_t^{x,c,\pi})_{ 0 \leq t \leq T}$  associated to the triple $( x,c, \pi)$  where $ x$ is the initial wealth, $\pi$ is the portfolio strategy and $c$ the consumption plan $c$, follows the dynamics given by:
\begin{eqnarray}
d X_t^{x,c,\pi}= \pi_tdS_t -c_tdt,\,\,\,\,\, X_0^{x,c,\pi}=x,
\end{eqnarray}
  The set of consumption-investment strategies  $(c, \pi)$  satisfying the   following no-bankruptcy condition is called
the admissible strategies set and  is denoted by $ \mathcal A (x)$ :
\begin{equation}
\P-a.s., \qquad  X_t^{x,c,\pi} \geq 0, \qquad \forall t \in [0,T].
\end{equation}
\\For this model, one assumes that :  \begin{itemize} \item The   appreciation rates $ ( \mu^{i},
i=1,\dots, d+p )$ are bounded predictable processes.\item The processes $(\varphi^{i,j}
 , 1\le i\le d+p, 1\le j\le  d)$  are bounded and  predictable and satisfy $\varphi^{i,j}_t>-1$
a.s. \item The
processes $(\sigma^{i,k} , 1\le i\le d+p, 1\le j\le  d)  $ are bounded and  predictable  \item The
matrix $\Sigma  $ is invertible. This condition  ensures that the market is arbitrage free.
\end{itemize}
Then,  the  pair consumption-terminal  wealth   satisfies the budget constraint
$$\Es^{\widetilde
{\mathbb{P}}}\big[\int_0^T c_t dt+ X_T^{ x,c,\pi}
  \big]\le x
$$
where $ \widetilde{\mathbb{P}}$ is the unique  equivalent   martingale measure 
with density  $ \wt Z $ given by \eqref{risk}, where  $(\theta_t,\gamma_t):= \Sigma_t ^{-1}\mu_t$ and $e^{-z^i}-1=\gamma_i$.
\end{exemple}


\subsection{Properties of the value process}
\noindent In this section, we  derive necessary conditions   satisfied by the optimal  control parameters. We start by showing regularity properties for the value process $V^{x,c,\psi}$ with respect to $(c, \psi)$.

\begin{proposition}\label{continuous} Define the map
 $G:\mathcal{A}(x)\longrightarrow D^{\exp}_0$
  as $G(c,\psi)=V$,  where
$(V,M^{V,c},v)$ is the solution of the BSDEJ \eqref{eqback}
associated with
$(U(c),\bar U(\psi))$.  Then \\[0.2cm]
(i) $G$ is concave, i.e., for all $\theta \in (0,1)$ and
$(c^1,\psi^1),(c^2,\psi^2) \in {\cal A}(x)$:
$$G\left(\theta c^1+(1-\theta)c^2,\theta \psi^1+(1-\theta)\psi^2\right)\ge \theta G(c^1,\psi^1)+(1-\theta)G(c^2,\psi^2).$$
(ii)  Let $G_0(c,\psi)$ be the value of
$G(c,\psi)$ at time 0, i.e., $G_0(c,\psi)=V_0$. If $(c^n,\psi^n)\in {\cal
A}(x)$ converges decreasingly to $(c,\psi)\in {\cal A}(x)$, then
$G_0( c^n,\psi^n)$ converges decreasingly to $G_0(c,\psi)$.
Moreover $G_0$ is upper continuous with respect to the control
parameters.
\end{proposition}
\noindent \textbf{Proof:} Let
$(V^k,M^{k,c},\vh^k )$ be the solution of the BSDE \eqref{eqback}
associated with $(U(c^k),\bar U(\psi^k))$ for $k = 1, 2$. For any
$\theta \in (0,1)$, let $ ( \wt {V}^\theta, \wt {M}^{
\theta,c}, \wt \vh^{\theta})$  be the solution of
\eqref{eqback} associated with $ (U(\theta c^1+(1-\theta)c^2),\bar
U(\theta \psi^1+(1-\theta)\psi^2))$ and $ ({V}^\theta,{M}^{
\theta,c},{\vh^{\theta}})$  be the solution of \eqref{eqback}
associated with   $\theta U(c^1)+(1-\theta)U(c^2),\theta \bar
U(\psi^1) +(1-\theta)\bar U(\psi^2))$  and   set
$\overline{V}^\theta :=\theta V^1+(1-\theta)V^2$. Then, by using
both the concavity  properties of ($U$, $\bar U$) and Theorem
\ref{comparaisonth},  we get $\widetilde V^\theta\ge V^\theta$.
Moreover, as consequence of Theorem \ref{concave}, we obtain
$V^\theta\ge
\overline{V}^\theta$, which gives the assertion $(i)$. \\
Let us now consider $(c^n,\psi^n)$ a decreasing  sequence of
control parameters in $\mathcal{A} (x)$ such that $c^n_t\longrightarrow c_t, \forall t$ a.s  and $\psi^n\longrightarrow
\psi$ a.s; then, by using inequality (\ref{comparison}), and the
fact that the functions $U$ and $\bar U$ are non-decreasing, we
get
\begin{equation}\label{exit}
{|V^{c^n,\psi^n}_0-V^{c,\psi}_0|} \le
\mathbb{E}^{\mathbb{Q}^*}\left[\int_0^T({U(c^n_s)-U(c_s)}) ds
+{(\bar U(\psi^n)-\bar U(\psi))}\right]
\end{equation}
where $\mathbb{Q}^*$ is the optimal measure associated with
$(U(c),\bar U(\psi))$. Thus,  by using the  monotone convergence theorem and the a priori estimate \eqref{priori},
$V^{c^n,\psi^n}$ converges decreasingly to $V^{c,\psi}$.
Let  $(c^n,\psi^n)\in {\cal
A}(x)$ be  a sequence of control parameters  such that
$c^n\longrightarrow c$ a.s. and $\psi^n\longrightarrow \psi$ a.s
where $(c,\psi)\in {\cal A}(x)$ and   denote $\widetilde
c^n=\sup_{m\ge n}c^m$ ,$\widetilde \psi^n=\sup_{m\ge n}\psi^m$.
Then, $\widetilde c^n\longrightarrow c$ a.s. decreasingly and
$\widetilde \psi^n\longrightarrow \psi$ a.s decreasingly. It
follows that   $V^{\widetilde c^n,\widetilde \psi^n}_0$ converges
to $V^{c,\psi}$ decreasingly and therefore:
$$\lim_{n}\sup V^{c^n,\psi^n}_0\le \lim_{n}V^{\widetilde c^n,\widetilde \psi^n}_0=V^{c,\psi}_0\,.$$
Hence, $G_0$ is upper semicontinuous with respect to the
control parameters.
$\hfill\Box$ 
\begin{definition} The pairs  $(c^1,\psi^1),(c^2,\psi^2) \in {\cal A}(x)$  are comparable if either
$(c^1,\psi^1 )\geq (c^2 ,\psi^2 )$ or $(c^1,\psi^1 )\leq (c^2
,\psi^2 )$ with the order introduced in Definition \ref{ordre}.
\end{definition}

\begin{proposition}\label{differ} Assume  that  Assumption A. \ref{utility} holds and let
$(c^1,\psi^1),(c^2 ,\psi^2) $ be two    comparable plans in $
  {\cal A}(x)$. Then the function $\Psi $, defined on $(0,1)$ and valued in $D_0^{\exp}$,  $$\Psi (\epsilon)  =G(c^1+\epsilon(c^2-c^1),\psi^1+\epsilon(\psi^2-\psi^1)) $$
is right-continuous at $0$.
\end{proposition}

\noindent \textbf{Proof:} Assume first that
$(c^1,\psi^1)\leq
 (c^2 , \psi^2)$. Let, for $\epsilon \in ]0,1[$,
$V^\epsilon=G(c^1+\epsilon(c^2-c^1),\psi^1+\epsilon(\psi^2-\psi^1))$
and $V =G(c^1,\psi^1)$.  From Proposition \ref{estimation} and the
obvious inequalities $U(c^1+\epsilon(c^2-c^1))\ge U(c^1)$ and
$\bar U(\psi^1+\epsilon(\psi^2-\psi^1))\ge \bar U(\psi^1)$, we
obtain
\begin{equation*}
\begin{split}
\mathbb{E}^{\mathbb{Q}^{*,2}}\Big(\sup_{0\le t\le
T}{|V _t-V^\epsilon_t|}^2 \Big)&\le C
\mathbb{E}^{\mathbb{Q}^{*,2}}\Big[{|\bar
U(\psi^1+\epsilon(\psi^2-\psi^1))-\bar U(\psi^1)|}^2\\& +\int_0^T
{|U(c^1_s+\epsilon(c^2_s-c^1_s))-U(c^1_s)|}^2 ds \Big].
\end{split}
\end{equation*}
\noindent Using now the   concavity properties of $U$ and $\bar U$, we obtain
 \begin{eqnarray*}0&\le& U(c^1_t+\epsilon(c^2_t-c^1_t))-U(c^1_t)\le \epsilon U'(c^1_t)(c^2_t-c^1_t) \\ 0&\le&\bar U(\psi^1+\epsilon(\psi^2-\psi^1))-\bar U(\psi^1)\le\epsilon\bar U'(\psi^1)(\psi^2-\psi^1).\end{eqnarray*} Thus, we have $$\mathbb{E}^{\mathbb{Q}^{*,2}}\Big(\sup_{0\le t\le
T}{\Big|{V _t-V^\epsilon_t\over \epsilon}\Big|^2}\Big)\le C
\mathbb{E}^{\mathbb{Q}^{*,2}}\Big[  { (\bar
 U'(\psi^1)})^2{(\psi^2-\psi^1)}^2+\int_0^T
({U'(c^1_s)})^2{(c^2_s-c^1_s)}^2 ds\Big].
$$

\noindent   Note that since  $(c^1,\psi^1)$ and $ (c^2 ,\psi^2) $ belong to  $
  {\cal A}(x)$, the right side term  of the last estimate is finite. Assume now that $(c^1,\psi^1)\geq (c^2, \psi^2)$. Then,
using the fact that $G$ is concave with respect to the control
parameters, one has $V^\epsilon\ge (1-\epsilon)V^{ 1}+\epsilon V^{ 2}$ where the $V^k$
are associated with $(c^k,\psi^k)$, hence $  {V^\epsilon-V^{1}\over \epsilon}\ge V^{2}-V^{1}.$  Moreover, since $c^1+\epsilon(c^2-c^1)\le c^1$ and
$\psi^1+\epsilon(\psi^2-\psi^1)\le \psi^1$, we have by Theorem
\ref{comparaisonth} that $0\ge {V^\epsilon-V^{1}\over \epsilon} .$
Therefore $\Big\vert {V^\epsilon_t-V^{1}_t\over
\epsilon}\Big\vert \le \vert  V^{2}_t-V^{1}_t\vert, t\in
[0,T].$ Using now  Proposition \ref{estimation}, we get
$$\label{pri}
 \mathbb{E}^{\mathbb{Q}^{*,1}}\Big(\sup_{0\le t\le
T}{\Big|\,{V^1_t-V^\epsilon_t\over \epsilon}\Big|}^2\Big)\le  c
\mathbb{E}^{\mathbb{Q}^{*,1}}\Big[ \, { (\bar
U'(\psi^2)})^2{(\psi^2-\psi^1)}^2+\int_0^T
({U'(c^2_s)})^2{(c^2_s-c^1_s)}^2 ds\Big]. $$
\noindent Finally, we conclude there exists a constant $C>0$ such that:
$\mathbb{E}^{\mathbb{Q}^*}\Big[\sup_{0\le t\le T}{\Big|{V^1_t-V^\epsilon_t\over
\epsilon}\Big|}^2\Big]\le C$,
 where $\mathbb{Q}^*=\mathbb{Q}^{*,1}$ if $(c^1,\psi^1 )\geq (c^2, \psi^2 )$
 and  $\mathbb{Q}^*=\mathbb{Q}^{*,2}$ if $(c^1,\psi^1 )\leq (c^2 , \psi^2 )$, then by Kolmogorov's criteria,
 we deduce that $\Psi$ is right-continuous at  $0$. $\hfill \Box$ \\\\
We give now a differentiability  result for our BSDE \eqref{eqback} .   We note that Imkeller et al. \cite{Imkeller10} have showed a  differentiability result  for a quadratic BSDE's driven by a continuous martingale, but the paper does not cover our case, so the proof of the result  will given in the Appendix.
\begin{theoreme}\label{difdym}
Let $(c^1,\psi^1)$ and $(c^2,\psi^2)$ be two comparable
plans in  ${\cal A}(x)$. Let  $
(V^\epsilon,M^{
\epsilon,c},v^\epsilon) $ be the solution of
 \eqref{eqback} associated with $(U(c^1+\epsilon(c^2-c^1)),\bar
U(\psi^1+\epsilon(\psi^2-\psi^1)))$ and $
(V^1,M^{ 1,c},v^1) $ the solution of
 \eqref{eqback} associated with $(U(c^1),\bar U(\psi^1))$. Then, $V^\epsilon$ is
right-differentiable with respect to $\epsilon$ at $0$. Moreover,
if we denote by $\partial_\epsilon V:= \displaystyle \lim_{\epsilon\rightarrow
0}{V^\epsilon-V^1\over \epsilon}$, then
 there exists $\partial_\epsilon \widetilde
M^{V,c},\partial_\epsilon v \in \mathbf{L}^2(\Q^{1,*})\times
\mathcal{L}^2(\widetilde\lambda ,\Q^{1,*})$
 such that the triple $(\partial_\epsilon  V,\partial_\epsilon  \widetilde M^{V,c},\partial_\epsilon v)$
 is the solution of the following BSDEJ:
\begin{equation}\label{derive}
\left\lbrace
\begin{split}
&d\partial_\epsilon V_t=\left(\delta_t\partial_\epsilon V_t
-U'(c^1_t)(c^2_t-c^1_t)\right)dt +d\partial_\epsilon \widetilde M^{V,c}_t +\sum_{i=1}^d \partial_\epsilon \vh^i_t d\widetilde N^i_t, \quad \Q^{*,1}\textbf{-}a.s.\\
&\partial_\epsilon V_T=\bar U'(\psi^1)(\psi^2-\psi^1),
\end{split}
\right.
\end{equation}
where $ \widetilde{\lambda}^i := \lambda^i e^{-v^{1,i}}$ and
$\widetilde N^i:=N^i-\int_0^. (e^{-v^{1,i}_t}-1)\lambda^i_tdt$ is
a $\Q^{1,*}$-martingale.\\\\\noindent Moreover, we obtain for all $t\le T$:
\begin{equation}\label{deriv_back}
\partial_\epsilon
V_t=\mathbb{E}^\mathbb{P}\left[{Z^{\Q^{*,1}}_T\over
Z^{\Q^{*,1}}_t}{S^\delta_T\over S^\delta_t}\bar
U'(\psi^1)(\psi^2-\psi^1)+\int_t^T {Z^{\Q^{*,1}}_s\over
Z^{\Q^{*,1}}_t}{S^\delta_s\over
S^\delta_t}U'(c^1_s)(c^2_s-c^1_s)ds \Big{\vert} {\cal G}_t\right].
\end{equation}
\end{theoreme}
\subsection{The Dynamic maximum principle}
\label{maximumprinciple}
\noindent We recall that we are interested in the following optimization problem: we associate with a pair $(c,\psi) \in
{\cal{A}} (x)$  the quantity $ X^{c,\psi}_0 =\Es^{\tilde \P} \left [\int_0^T c_sds+ \psi
\right]$  and we study
\begin{equation}\label{fopt}
u(x)=\sup_{X^{c,\psi}_0\le x}V^{(c,\psi)}_0 \,.\end{equation} Here
$V^{(c,\psi)}_0=V_0$, where $(V,M^{V,c},\vh)$ is the solution of
the BSDE (\ref{eqback}) associated with $(U(c),\bar U(\psi))$.
Note that if we are in the setting of  Example \ref{example2}, our problem  correspond to a maximization  of recursive utility function  over consumption-investment strategy where  $X_0$ is the initial wealth associated with the consumption $c$
and terminal wealth $\psi$.

\begin{proposition}
 There exists an optimal pair $(c^0,\psi^0)$ which solves (\ref{fopt}).
 \end{proposition}

\noindent
\textbf{Proof}: The uniqueness is a consequence of the strictly concavity property of $V_0$. We shall prove the existence by using Koml\`os theorem.
\\[0.2cm] \noindent \textbf{First step}:
Let us first prove that  $\sup_{(c,\psi)\in {\cal
A}(x)}V^{c,\psi}_0 <+\infty$. Because $\P\in {\cal Q}^e_f$, we
have:
$$\sup_{(c,\psi)\in {\cal A}(x)}V^{c,\psi}_0\le \sup_{(c,\psi)\in {\cal A}(x)}\Es^\P\Big[\bar U(\psi)+\int_0^T U(c_s)ds \Big] =:\widetilde u(x)$$
Using the elasticity assumption on $U$ and $\bar U$, we can find
$\gamma\in (0,1)$ and  $x_0\in \R$ such that, for any $\theta >1$,
one has:
\begin{equation*}
U(\theta x)<\theta^\gamma U(x),\qquad
 \bar U(\theta x)<\theta^\gamma \bar U(x) \qquad  \forall x>x_0,
\end{equation*}
hence,   for any $x>x_0$:
$$\widetilde u(\theta x)=\Es^\P\Big[\bar U(\theta{\psi^{\theta x}\over \theta})+\int_0^T U(\theta {c^{\theta x}_s\over \theta})ds \Big]\le \theta^\gamma \widetilde u(x).$$
Then, $AE(\widetilde u)<1$, which permits to conclude that, for any $x>0$ , $\widetilde u(x)<+\infty$
(see \cite{KS99} and  \cite{M09} chap. 3, Lemma 3).\\
\\
\noindent
\textbf{Second step:} Let $(c^n,\psi^n)\in {\cal A}(x)$ be a maximizing sequence   such that:
$$ \nearrow  \lim_{n \rightarrow  +\infty} V^{c^n,\psi^n}_0= \sup_{(c,\psi)\in {\cal A}(x)} V^{c,\psi}_0<+\infty,$$
where the  RHS    is finite thanks to step 1. Using Koml\`os
criterion, we can find a convex combination $(\bar c^n,\bar
\psi^n)\in conv\Big((c^n,\psi^n),(c^{n+1},\psi^{n+1}),\cdots\Big)$
which converges $\mathbb{P}\textbf{-}a.s$. We denote  by
$(c^*,\psi^*)$ this limit,  which belongs to  ${\cal A}(x)$ since
it is a  closed convex set.  Moreover, there exists $N_n\ge n$ and
a positive sequence ${(\theta^m)}_{m\in \mathbb{N}}$ satisfying
$\sum_{m=n}^{N_{n}} \theta^{m} =1$ such that $(\bar c_n, \bar
\psi_n) =(\sum_{m=n}^{N_n} \theta^m c^m,\sum_{m=n}^{N_n} \theta^m
\psi^m)$. Therefore, by using the concavity and the increasing
properties of the functional $V_0$ which respect to the control
plan we get:
$$V^{\bar c^n,\bar \psi^n}_0=V^{\sum_{m=n}^{N_n}\theta^m c^m,\sum_{m=n}^{N_n}\theta^m \psi^m}_0\ge \sum_{m=n}^{N_n}\theta^m V^{c^m,\psi^m}_0\ge V^{c^n,\psi^n}_0.$$
Moreover,  using the upper semi-continuous property of the
functional $V_0$ which respect to  the control  plan we get:
$$\sup_{(c,\psi)\in \mathcal{A}(x)}V^{c,\psi}_0
=\lim\sup_{n}V^{c^n,\psi^n}_0\le \lim\sup_n V^{\bar c^n,\bar \psi^n}_0=V^{c^*,\psi^*}_0.$$
$\hfill\Box$\\
\noindent
 In order to characterize the optimal solution, we recall the
classical convex analysis result.
\begin{proposition}
There exists a constant $\nu^*>0$ such that :
\begin{equation}\label{seconopt}
u(x)=\sup_{(c,\psi)}\Big\{V^{c,\psi}_0+\nu^*\left(x-X^{c,\psi}_0\right)\Big\}
\end{equation}
\noindent and if the maximum is attained in (\ref{fopt})  by $(c^*,\psi^*)$,
then it is attained in $(\ref{seconopt}) $ by $(c^*,\psi^*)$ with
$X^{c^*,\psi^* }_0=x$. Conversely, if there exists $\nu^0>0$ and
$(c^0,\psi^0)$ such that the maximum is attained in
$$\label{thirdopt}
\sup_{(c,\psi)}\Big\{V^{c,\psi}_0+\nu^0\left(x-X^{c,\psi}_0\right)\Big\}
$$
 with $X^{c^0,\psi^0}_0=x$, then the maximum is attained in (\ref{seconopt}) by $(c^0,\psi^0)$.
\end{proposition}
 Let  $\nu>0$ be fixed and  $L$ be the map   given by $L(c,\psi)=V^{c,\psi}_0-\nu X^{c,\psi}_0$. We now study  the following optimization problem:
\begin{equation}\label{Lan}
\sup_{(c,\psi)}L(c,\psi)\,.
\end{equation}

 \begin{proposition}\label{optidym}
The optimal   plan $(c^0,\psi^0)$ which solves
(\ref{Lan}) satisfies the following (implicit) equations:

\begin{equation}\label{consopt}
 U'(c^0_t)={Z^{\widetilde \P}_t\over Z^{0}_t}{\nu\over
S^\delta_t}\quad  dt\otimes d\P  \hbox{ a.s },\qquad \bar
U'(\psi^0)={Z^{\widetilde \P}_T\over Z^{ 0}_T }{\nu\over
S^\delta_T},\, d\P  \hbox{ a.s }
\end{equation}
where  $Z^0$ is the Radon-Nikodym density of  the probability
measure $\Q^0$   associated with the optimal plan $(c^0,\psi^0)$.
\end{proposition}

\noindent \textbf{Proof:}  Consider the optimal  plan
$(c^0,\psi^0)$ which solves (\ref{Lan}) and another
plan $(c,\psi)$. For  $\epsilon\in (0,1)$, one has $L(c^0+\epsilon(c-c^0),\psi^0+\epsilon(\psi-\psi^0))\le L(c^0,\psi^0)$, then
\begin{equation}\label{difLan}
{1\over \epsilon}\left({V^{c^0+\epsilon(c-c^0),\psi^0+\epsilon(\psi-\psi^0)}_0-V^{c^0,\psi^0}_0}\right)-\nu{1\over \epsilon}\left({X^{c^0+\epsilon(c-c^0),\psi^0+\epsilon(\psi-\psi^0}_0-X^{c^0,\psi^0}_0}\right)\le 0
\end{equation}From the definition, we obtain that
$$\partial_\epsilon X^{c^0,\psi^0}_0: =\lim_{\epsilon\rightarrow
0}{1\over
\epsilon}{(X^{c^0+\epsilon(c-c^0),\psi^0+\epsilon(c-c^0)}_0-X^{c^0,\psi^0}_0)}=\Es^{\widetilde
\P}\left[\int_0^T (c_s-c^0_s)ds +(\psi-\psi^0) \right]\,.$$ Taking
the limit when $\epsilon $ goes to 0   in (\ref{difLan}), we
obtain:
\begin{equation}\label{difsuite}
\partial_\epsilon V^{c^0,\psi^0}_0-\nu \partial_\epsilon X^{c^0,\psi^0}_0\le 0
\end{equation}
where $ \partial_\epsilon V^{c^0,\psi^0)} $ exists and   is
given explicitly by  Theorem \ref{difdym}. Note that:
\begin{eqnarray*}
\partial_\epsilon V^{c^0,\psi^0}_0-\nu \partial_\epsilon
X^{c^0,\psi^0}_0&=& \Es^\P \left[ S^\delta_T Z^{0}_T \bar
U'(\psi^0)(\psi-\psi^0)
+\int_0^T S^\delta_s Z^{0}_s U'(c^0_s)(c_s-c^0_s)ds\right]\nonumber\\
&-&\nu \Es^\P\left[  Z^{\widetilde \P}_T(\psi-\psi^0)+\int_0^T
Z^{\widetilde \P}_s(c_s-c^0_s)ds \right]\label{d1}
\end{eqnarray*}
It follows from the equality (\ref{difsuite}) that  $$\label{d2}
\Es^\P \left[ \left(S^\delta_T Z^{0}_T  \bar U'(\psi^0)-\nu \mb{
Z^{\widetilde \P}_T}\right)(\psi-\psi^0)+\int_0^T \left(S^\delta_s
Z^{0}_s U'(c^0_s)-\nu Z^{\widetilde
\P}_s\right)(c_s-c^0_s)ds\right ]\le 0
$$
The \mb{end} of the proof is the same as  in El Karoui et al.
\cite{ElkPQ01} (proof of Theorem 4.2, p. 677). In particular,  for
any $\psi$, $\Es^\P \left[  \left(S^\delta_T Z^{0}_T  \bar
U'(\psi^0)-\nu Z^{\widetilde \P}_T\right)(\psi-\psi^0) \right]\le
0$, hence
$$ S^\delta_TZ^{0}_T\bar U'(\psi^0)-\nu Z^{\widetilde \P}_T = 0 \qquad a.s$$
We find the optimal $c$ with similar arguments.
$\hfill\Box$

\begin{theoreme}\label{finalopt} Let $I$ and $\bar I$ be the inverse of the functions $U'$ and $\bar U'$. The
 optimal plan $(c^0,\psi^0)$ which solve the   problem (\ref{seconopt}) is given by:
\begin{equation*}\label{d3}
c^0_t=I\left({\nu^0\over  S^\delta_t}{Z^{\widetilde \P}_t\over
Z^{0}_t}\right) \qquad dt\otimes d\P \hbox{ a.s },\qquad
\psi^0=\bar I\left({\nu^0\over S^\delta_T}{Z^{\widetilde
\P}_T\over Z^{0}_T}\right)\qquad \P  \hbox{ a.s. }.
\end{equation*}
where $\nu^0>0$ satisfies:
\begin{equation*}\label{d3a}
\mathbb{E}^{\widetilde \P}\left[\int_0^T I\left( {\nu^0\over
S^\delta_t} {Z^{\widetilde \P}_t\over Z^{0}_t}\right)dt+\bar
I\left( {\nu^0\over S^\delta_T} {Z^{\widetilde \P}_T\over
Z^{0}_T}\right ) \right]=x.
\end{equation*}
\end{theoreme}
 \noindent \textbf{Proof:}
Define the map:  $f:(0,+\infty)\rightarrow (0,+\infty)$ as
$$f(\nu)= \mathbb{E}^{\widetilde \P}\left[\int_0^T I\left(
{\nu\over S^\delta_t} {Z^{\widetilde \P}_t\over
Z^{0}_t}\right)dt+\bar I\left( {\nu\over S^\delta_T}
{Z^{\widetilde\P}_T\over Z^{0}_T}\right ) \right] \,.$$ Then,
using assumption A.\ref{utility}, $f$ is monotone and satisfies
$\lim_{\nu\rightarrow 0}f(\nu)=+\infty$ and $\lim_{\nu\rightarrow
+\infty}f(\nu)=0$. For any  initial wealth $x\in (0,+\infty)$,
there exists a unique $\nu^0$ such that $f(\nu^0)=x$. \\\noindent
Let $(c,\psi)\in {\cal A}(x)$ and $(V^{c,\psi},M^{V,c},\vh)$
$\big(\hbox{resp.}\, (V^{c^0,\psi^)},M^{V^0,c},\vh^0) \big)$ the
solution of the BSDEJ (\ref{eqback}) associated with $(U(c^0),\bar
U(\psi^0))$ $\big(\hbox{resp. }\,(U(c),\bar U(\psi))\big)$ then
from the inequality (\ref{comparison}) (see the comparison
theorem), we get:
\begin{eqnarray*}\label{compfin}
V^{c,\psi}_0-V^{(c^0,\psi^0)}_0 & \le &
\mathbb{E}^{\Q^0}\left[S^\delta_T \left(\bar U(\psi)-\bar
U(\psi^0)\right)+\int_0^T S^\delta_s
\left(U(c_s)-U(c^0_s)\right)ds\right]\\& \le&
\mathbb{E}^{\Q^0}\left[ S^\delta_T \bar
U'(\psi^0)(\psi-\psi^0)+\int_0^T S^\delta_s
U'(c^0_s)(c_s-c^0_s)ds\right].
\end{eqnarray*}
It follows that:
\begin{eqnarray*}\label{compfina}
V^{c,\psi}_0-V^{(c^0,\psi^0)}_0 & \le& \nu^0
\mathbb{E}^{\Q^0}\left( {Z^{\widetilde \P}_T\over
Z^{0}_T}(\psi-\psi^0)+\int_0^T
  {Z^{\widetilde \P}_s\over Z^{0}_s}(c_s-c^0_s)ds\right)\\
& \le &\nu^0\left(\mathbb{E}^{\widetilde \P}\left( \psi+\int_0^T  c_s
ds\right)-\mathbb{E}^{\widetilde \P}\left( \psi^0+\int_0^T  c^0_s
ds\right)\right)
\end{eqnarray*}
Since $(c,\psi) \in {\cal A}(x)$, then $\mathbb{E}^{\widetilde
\P}\left[  \psi+\int_0^T c_s ds\right]\le x$. Using that
$\mathbb{E}^{\widetilde \P}\left[\psi^0+\int_0^T  c^0_s
ds\right]=x$,  we conclude:
$$ V^{c,\psi}_0\le V^{c^0,\psi^0}_0\,.$$
$\hfill\Box$

\section{ Logarithm  Case}
\noindent In this section, we   assume that the process $\delta$ is
deterministic    and that $U(x)=\ln (x)$ and $\bar U(x)=0$, hence   $I(x)= \displaystyle {1\over x}$ for all $x\in (0,+\infty)$.
We introduce, as   in Theorem \ref{finalopt}, the optimal process
$c^*_t= \displaystyle I\Big({\nu\over   S^\delta_t}{ \widetilde Z _t\over Z^
* _t}\Big)={  S^\delta_t \over\nu }{
 Z^* _t \over \widetilde Z _t} $. Recall  that the Radon-Nikodym
density    $\wt Z$, and the Radon-Nikodym density of the optimal
probability measure $ Z^*$ (given   in \eqref{densitel}) satisfy
\beq d\wt Z_t&=& \wt Z_{t-}
(\theta_tdM^c_t+\sum_{i=1}^n (e^{-z^i_t}-1)dN^i_t),\, \wt Z_0=1\\
d  Z^*_t&=&   Z^* _{t-} (-dM^{Y,c}_t+\sum_{i=1}^n
(e^{-y^i_t}-1)dN^i_t),\, Z^*_0=1.\eeq
For any deterministic function $\alpha$ such that $\alpha(T)=0$,
$V$ admits a decomposition as
$$V_t= \alpha (t) \ln (c^*_t)+ \beta_t$$ where $\beta$ is a
process such that $\beta _T=0$.  Our goal is to characterize the
process $\beta$. As in \cite{BMS07th}, we introduce $J_t
=\frac{1}{1+\alpha(t)} \beta_t$ in order to obtain a simple BSDEJ.
Note that, even if $Z^*$  is implicit (the coefficients depend on
the solution $c^*$), the BSDEJ for  $J$ is explicitly determined in
terms of the given parameters $\lambda^i$ and of the given
probability $\widetilde \P$.
\begin{proposition}
The value function $V$ has the form $V_t= \alpha (t) \ln (c^*_t)+ (1+\alpha (t)) J_t$ where
$\alpha(t) = -\int_t^T e^{\int_t^s \delta (u)du}ds$  and  $(J,\bar M^{J,c},j)$ is the unique
solution of the following BSDEJ:
\begin{equation}\label{forward_J}
\left\lbrace
\begin{split}
&dJ_t= \Big(\, (1+\delta
(t))(1+k(t))  J_t -k(t)\delta(t) \Big)dt +  d \bar M^{J,c}_t
+\frac 12
  d\cro{ \bar M^{J,c}}_t+ \frac 12 k(t)(1+k(t))\theta_t^2 d\langle M
  ^c\rangle_t
 \nonumber \\&\hspace{1cm}+\sum_{i=1}^d j ^i_t
d\bar N^i_t +\sum_{i=1}^d
  \Big(g(j^i_t) \bar \lambda ^i_t +  \Big (k(t)(e^{-z^i_t}-1)+ e^{k(t)z^i_t}-1  \Big)\lambda
^i_t  \Big)   dt\\& J_T= 0
\end{split}
\right.
\end{equation}
where $k(t)=-\frac{\alpha(t)}{ 1+\alpha (t) }$.
  Here,  the processes  $\bar M^{J,c} $ and $d\bar N_ t^i= dH^i_t -\bar \lambda ^i_t dt$ are $\bar \P$-martingales
where  $ d\bar\P \vert_{{\cal G}_t}= \bar Z_t d\P \vert_{{\cal
G}_t}$, $\bar \lambda ^i_t = e^{k(t) z^i_t} \lambda ^i_t$ and
 \beq  d \bar Z _t&=&  - \bar Z _{t^-} \left(   k(t)
 \theta _tdM^c_t - \sum_ {i=1}^d(e^{k(t) z^i_t}-1) dN^i_t  \right)\eeq
   \end{proposition}
\noindent Note that, in a complete market, one obtains a forward backward system for the pair $J$-optimal wealth.
\\\\
\noindent \textbf{Proof:}  Using the fact that $V$ satisfies the
BSDE (\ref{eqback}) and the assumed form of $V$ in terms of
$(\alpha,\beta)$, one obtains $dV_t=(\delta(t)V_t - \ln (c^*_t))\,dt-d(\ln Z^*_t)=\alpha (t) d(\ln c^*_t)+ (\ln c^*_t) \alpha^\prime (t)dt
+d\beta_t$. Therefore \beqa d\beta_t&=& \delta(t)(V_t + \alpha
(t))  dt - (1+\alpha^\prime(t))\ln (c^*_t)dt + \alpha(t) d\ln \wt
Z_t +(\alpha(t)+1) d\ln Z^*_t\\&=&\big(
\left(\delta(t)\alpha(t)-1-\alpha^\prime(t)\right) \ln c^*(t)
+\delta (t)\beta_t + \alpha (t) \delta (t)  \big)dt+ \alpha(t)
d\ln \wt Z_t +(\alpha(t)+1) d\ln Z^*_t\eeqa We choose $\alpha$ so
that $ \delta(t)\alpha(t) =1+\alpha ^\prime (t)$. It follows that
$$d\beta_t= \delta(t)(\beta_t + \alpha (t) ) dt
  +\alpha(t) d\ln \wt Z_t
+(\alpha(t)+1) d\ln Z^*_t$$ After some obvious computations taking
into account the form of $\wt Z$ and $ Z^*$, one obtains \beqa
d\beta_t&=& \delta(t)(\beta_t + \alpha (t)) dt + \sum_{i=1}^d
\Big((\alpha(t)+1) (e^{-y^i_t}-1)-\alpha(t)
(e^{-z^i_t}-1)\Big)\lambda ^i_t dt\\&&+ \alpha (t)\theta_t
dM^c_t+(\alpha(t) +1)dM^{V,c}_t   -\frac 12 \Big(\alpha
(t)\theta^2_t d\cro{M^c}_t -(\alpha(t)+1) d\cro{M^{V,c}}_t\Big)
\\&&+\sum_{i=1}^d \Big((\alpha (t)+1)y ^i_t -\alpha (t) z^i_t \Big) dH^i_t\eeqa
We now define $J_t:= \frac{1}{1+\alpha(t)}\beta_t$ and set
$k(t)=-\frac{\alpha (t)}{1+\alpha(t)}$, then  we find the following dynamics: \beqa dJ_t&=&
\left(\frac{ 1+\delta (t)  }{1+\alpha(t) } J_t -\delta(t)k(t)
\right)dt + \sum_{i=1}^d \Big(  g(y^i_t )+k(t)g(z^i_t)
\Big)\lambda ^i_t dt\\&&+   dM^{V,c}_t-k(t)\theta_t dM^c_t  +\frac
12 \Big(k(t) \theta^2_t d\cro{M^c}_t + d\cro{M^{V,c}}_t\Big)
 +\sum_{i=1}^d \Big( y ^i_t +k(t) z^i_t
\Big) dN^i_t\eeqa We introduce the martingale $M^{J,c}$  as
 $dM^{J,c}_t:=dM^{Y,c}_t-k(t)\theta_t dM^c_t$.
It is easy to check that $$d\langle M^{J,c}\rangle_t = d\langle
M^{Y,c}\rangle_t -k ^2(t)\theta^2_t d\langle M^{ c}\rangle_t
 -2k(t)\theta _t d\langle
M^{J,c}, M^c\rangle_t$$  and we denote $j^i_t= y ^i_t +k(t)
z^i_t$.  Using the fact that,  due to the form of $g$, for any
$x,k,z,\lambda$,
$$xdN_t + \lambda (g(x-kz)+kg(z))dt=x(dN_t-(e^{kz}-1)\lambda dt) +\left(g(x)
 e^{kz}+    (e^{-z}-1)k+e^{kz}-1\right)\lambda dt$$ one
obtains
 \beqa dJ_t&=&
\Big((1+\delta (t))(1+k(t))  J_t -\delta(t)k(t) \Big)dt +
\sum_{i=1}^d \Big(   g(j^i_t )  e^{k(t)z^i_t}+ k(t)(e^{-z^i_t}-1)+
e^{k(t)z^i_t}-1  \Big)\lambda ^i_t dt\\&&+ dM^{J,c}_t +\frac 12
d\langle M^{J,c}\rangle _t+k(t) \theta_td\langle M^{J,c},
M^c\rangle_t+\frac 12 k(t)(k(t)+1)\theta^2_t d\langle M^c\rangle
_t
\\&&+\sum_{i=1}^d j_  t^i (dN^i_t- (e^{kz^i_t}-1)\lambda ^i_tdt)\eeqa We define $\bar \P$  as $d\bar \P =
\bar Z  d\P$, where $d \bar Z _t= - \bar Z _{t-} \left [  k(t) \theta_t dM^c_t -   \sum_ {i=1}^d(e^{k(t) z^i_t}-1) dN^i_t  \right]$. The processes   $\bar M^{J,c}
$ and $\bar N^i$ defined as $d\bar M_t^{J,c} =dM_t^{J,c} +
 k(t) \theta_t  d\langle M^{J,c},M^c\rangle_t$ and  $d\bar N^i_t= dN^i_t -(e^{k(t) z^i_t}-1) \lambda ^i _tdt = dH^i_t -\bar \lambda ^i_t dt$ are $\bar \P$
 martingales. The result follows.

%

\section{Appendix : Proof of Theorem \ref{difdym}}
 \noindent  Let $  (V^\epsilon,M^{ \epsilon,c},v^\epsilon)$   be the solution of   \eqref{eqback} associated with $  (U(c^1+\epsilon(c^2-c^1)),\bar
U(\psi^1+\epsilon(\psi^2-\psi^1)))$  and $(V^1,M^{1,c},v^1) $ be the solution of   \eqref{eqback} associated with $ \big(U(c^1),\bar
U(\psi^1) \big)$  and denote
\begin{equation}
\label{derivation}
\begin{split}
& \Delta_\epsilon V:={V^\epsilon-V^1\over \epsilon},\qquad \Delta_\epsilon M^{ c}:={M^{ \epsilon,c}-M^{ 1,c}\over \epsilon},\qquad  \Delta_\epsilon v^i:={v^{\epsilon,i}-v^{1,i}\over \epsilon},\\
&\Delta_\epsilon U := \frac{U(c^1+\epsilon(c^2-c^1)) - U(c^1)}{\epsilon}, \qquad \Delta_\epsilon \bar U_T := \frac{\bar U(\psi^1+\epsilon(\psi^2-\psi^1)) - \bar U(\psi^1)}{\epsilon}
\end{split}
\end{equation}
\noindent then,  $ (\Delta_\epsilon V, \Delta_\epsilon M^{ c}, \Delta_\epsilon v)$ satisfies the following   equation:
\begin{equation}
\label{BSDE:approximation}
\begin{split}
&\Delta_\epsilon V_t-\int_0^t(\delta_s \Delta_\epsilon
V_s-\Delta_\epsilon U_s)ds={1\over 2 \epsilon }(\langle
M^{\epsilon,c}\rangle_t-\langle M^{1,c}\rangle_t)\\&+{1\over
\epsilon}\sum_{i=1}^d\int_0^t
(g(v^{\epsilon,i}_s)-g(v^{1,i}_s))\lambda^i_s ds +\Delta_\epsilon
M^{c}_t +\sum_{i=1}^d\int_0^t \Delta_\epsilon v^i_s dN^i_s,
\end{split}
\end{equation}
\noindent with final condition $\Delta_\epsilon V_T = \Delta_\epsilon \bar U_T$. We start first to give the following  a priori estimates:
 \begin{lemma}\label{control_mart} Assume the same conditions as in Theorem \ref{difdym}.
Then,  there exists a constant $C>0$ such that: $\forall \,  i=1,\cdots, d, \forall \, p\in \N^*$, $\forall \, \epsilon >0$,
\begin{equation}
\label{priori:estimate}
 \Es^{\Q^{*,1}}\left[\sup_{0\leq t \leq T} \left| \Delta_\epsilon  V_t\right|^2 +
  \left\langle \Delta_\epsilon \widetilde M^{ c}\right\rangle_T +\sum_{i=1}^d \int_0^T \frac{\left\vert \Delta_\epsilon v_s^i
\right\vert^p}{p!} \widetilde \lambda^i_s ds\right] \le C,
\end{equation}
  where $\Delta_\epsilon \widetilde M^{ c}$ is the $\Q^{*,1}$
martingale part of the $\Q^{*,1}$ semimartingale $\Delta_\epsilon
M^{ c}$, and $\widetilde \lambda^i :=\lambda^i e^{-v^{1,i}}$ is the intensity process of the process $H^i$ under the probability measure $\Q^{*,1}$.
\end{lemma}
\noindent
\textbf{Proof:}
 Let $(c^1,\psi^1)$ and
$(c^2,\psi^2)$ be two comparable  plans.
  We introduce   the processes
\begin{eqnarray*}
 K^\epsilon_t&:=&\Es^\P\left[\exp\left(\int_0^T \big(\delta_s V^\epsilon_s- U(c^1_s+\epsilon(c^2_s-c^1_s) )\big)ds- \bar U(\psi^1
+\epsilon(\psi^2-\psi^1))\right)\Big\vert {\cal{G}}_t\right]\\\\
K^1_t&:=&\Es^\P\left[\exp\left(\int_0^T \big(\delta_s V^1_s-
U(c^1_s)\big)ds- \bar U(\psi^1) \right)\Big\vert
{\cal{G}}_t\right].
\end{eqnarray*}
Obviously, for all $t\in [0,T]$, one has:
\begin{eqnarray*}
V^\epsilon_t&=&-\ln(K^\epsilon_t)+\int_0^t(\delta_s V^\epsilon_s- U(c^1_s+\epsilon(c^2_s-c^1_s) )ds\\
V^1_t&=&-\ln(K^1_t)+\int_0^t(\delta_s V^1_s- U(c^1_s) )ds\,,
\end{eqnarray*}
hence,
\begin{equation}\label{utildif}
\frac{V^\epsilon_t-V^1_t}{
\epsilon} =-\ln\left[\left({K^\epsilon_t \over
K^1_t}\right)^{1\over \epsilon}\right]+\int_0^t\left[\delta_s
\,\frac 1 \epsilon \left(V^\epsilon_s-V^1_s\right) -\frac 1 \epsilon \left( U(c^1_s+\epsilon(c^2_s-c^1_s))-U(c^1_s)\right) \right]ds\,.
\end{equation}
For $ t\in [0,T] $,  we define   $  \widetilde{K^\epsilon_t}  :=
\displaystyle \frac{ K^\epsilon }{ K^1_t}  $ and $
\bar{K}^\epsilon_t= \left({\widetilde{K^\epsilon_t}}\right)^{1/
\epsilon}$. The processes $\bar K^\epsilon$ and ${  (\bar
K^\epsilon)^{-1}}$ are positive  semi-martingales which belong  to
$\mathbf{{L}}^p(\P)$ since:
\begin{equation*}
{\left(
K^\epsilon_t\right)^{-p}}=\exp\left[p\Delta_\epsilon V _t+\int_0^t p(
\Delta_\epsilon U(c^1_s)-\delta_s \Delta_\epsilon V _s)ds\right].
\end{equation*}
In the other hand, by using the dynamics of $K^{\epsilon}$ and $K^1$ under the probability measure $\mathbb{P}$:
 \beqa
   {d K^{\epsilon}_t} &=&{K^{\epsilon}_{t-}}\left[ - dM^{\epsilon, c}_t + \sum_{i=1}^d \big( e^{- v^{\epsilon,i}_t } - 1 \big) dN_{t}^i\right]\\
    {d K^{1}_t} &=&{K^{\epsilon}_{t-}}\left[ - dM^{1, c}_t + \sum_{i=1}^d \big( e^{- v^{1,i}_t } - 1 \big) dN_{t}^i\right]
  \eeqa
 and applying  integration by parts formula, we get the dynamics of  $  \widetilde{K^\epsilon}$ given by:
\begin{equation}\label{submartingale:K}
\begin{split}
d  \widetilde{K^\epsilon_t}   &=   \widetilde{K^\epsilon_{t-}} \big[ - d \big( M^{ \epsilon,c}_t-M^{ 1,c}_t  - \langle M^{ \epsilon,c}-M^{ 1,c}, M^{ 1,c}\rangle_t \big)
\\&+ \sum_{i =1}^d \big( e^{- (v^{\epsilon,i}_t - v^{1,i}_t )} - 1 \big) \big[ dH_t^i - e^{ - v^{1,i}_t} \lambda_t^i dt \big] \, \big]
\end{split}
\end{equation}
Clearly, $\widetilde{K^\epsilon} $ is  $\Q^{*,1}$-local martingale.
Then,  the processes $\bar{K}^\epsilon$ and ${( \bar{K}^\epsilon)^{-1}}$ are
positive $\Q^{*,1}$-submartingales.
We now split the study into two cases. \\[0.2cm]
\noindent \textit{First case}: $(c^1, \psi^1)\leq (c^2, \psi^2)$. Using the inequality
\eqref{comparison}, for all $t\in [0,T]$:
\begin{equation*}\label{deltaV}
\begin{split}
&|\Delta_\epsilon V _t|\le \Es^{\Q^{*,1}}\left[{S^\delta_T\over
S^\delta_t}\bar U'(\psi^1)(\psi^2-\psi^1)+\int_t^T{S^\delta_s\over
S^\delta_t} U'(c^1_s)(c^2_s-c^1_s)ds\Big\vert {\cal G}_t\right]\\
&\sup_{0\le t\le T}{\left(1\over  \bar
K^\epsilon_t\right)}^p\le\exp\left[
p(||\delta||_{\infty}+1)\sup_{0\le t\le T} |\Delta_\epsilon V
_t|+\int_0^T p U'(c^1_s)(c^2_s-c^1_s)ds\right].
\end{split}
\end{equation*}
 Setting $ \kappa=p(c +1)$ where $c$  is the constant given in  Assumption A2, we obtain:
 \begin{equation*}
 \begin{split}
 \sup_{0\le t\le
T}{\left(1\over  \bar K^\epsilon_t\right)}^p
 & \le \exp\Big( \kappa \sup_{0\le t\le
T}\Es^{Q^{*,1}}\big[\bar U'(\psi^1)(\psi^2-\psi^1) +\int_t^T
U'(c^1_s)(c^2_s-c^1_s)ds\Big\vert{\cal G}_t\big] \\
&  \qquad +\int_0^T p
U'(c^1_s)(c^2_s-c^1_s)ds\Big).\\
\end{split}
\end{equation*}
\noindent Using Jensen's inequality, we have:
\begin{eqnarray}
 \sup_{0\le t\le T}{\left(1\over  \bar
K^\epsilon_t\right)}^p&\le &\sup_{0\le t\le
T}\Es^{\Q^{*,1}}{\left[\exp\left(\bar U'
(\psi^1)(\psi^2-\psi^1)+\int_t^T
U'(c^1_s)(c^2_s-c^1_s)ds\right)\Big\vert{\cal G}_t\right]}^\kappa
\nonumber\\ && \;\times  \exp(\int_0^T p U'(c^1_s)(c^2_s-c^1_s)ds).
\label{barK}
\end{eqnarray} Thanks to the assumption $(c^i,\psi^i) \in {\cal A}(x)$, we conclude that $\sup_{0\le t\le T}{\left(1\over
\bar{K}^\epsilon_t\right)} \in \mathbf{{L}}^p(\P)$.

\noindent \textit{Second case:} $(c^2, \psi^2)\geq (c^2, \psi^2)$.
Then, using concavity property, we obtain for all $t\in [0,T]$:
\begin{equation*}
\Big\vert {V^\epsilon_t-V^1_t\over \epsilon}\Big\vert\le |V^1_t-V^2_t|,\quad |\Delta_\epsilon U(c^1_t)|\le U'(c^2_t)(c^1_t-c^2_t)
\end{equation*}
Now,   using the same arguments as in the  first step,  we get that:
\begin{equation*}
\begin{split}
\sup_{0\le t\le T}{\left(1\over  \bar K^\epsilon_t\right)}^p\le &\sup_{0\le t\le T}\Es^{Q^{*,2}}{\left[\exp\left(\bar U' (\psi^2)(\psi^1-\psi^2)+\int_t^T  U'(c^2_s)(c^1_s-c^2_s)ds\right)\Big\vert{\cal G}_t\right]}^\kappa\\ & \; \times \exp(\int_0^T p U'(c^2_s)(c^1_s-c^2_s)ds).\\
\end{split}
\end{equation*}
 \noindent We use the same arguments to prove  $\sup_{0\le t\le
T}{|\bar K^\epsilon_t|} \in \mathbf{{L}}^p(\P)$. From the
representation theorem, there exists two continuous martingales
$\bar M^{\epsilon,c},\widetilde M^{\epsilon,c}$ and $d$
predictable processes $\bar k^{\epsilon,i},\widetilde k^{\epsilon,i}$ such
that:
\begin{equation*}
\begin{split}
&\bar K^\epsilon_t=K^\epsilon_0+\bar M^{\epsilon,c}_t+ \sum_{i=1}^d\int_0^t \bar k^{\epsilon,i}_s dN^i_s\\
&{1\over \bar K^\epsilon_t}={1\over K^\epsilon_0}+\widetilde
M^{\epsilon,c}_t+\sum_{i=1}^d\int_0^t \widetilde k^{\epsilon,i}_s dN^i_s\,.
\end{split}
\end{equation*}
These  processes being positive $\Q^{*,1}$-submartingales, using
(\ref{barK}) there exists two constants $C_K$ and $\widetilde C_K$
such that:
\begin{equation}\label{barsaut}
\begin{split}
&\Es^{\Q^{*,1}}\left[\int_0^T \bar {(k^{\epsilon,i}_s)}^2 \widetilde \lambda^i_s ds \right]\le \Es^{\Q^{*,1}}[{(\bar K^\epsilon_T)}^2]\le C_K\\
&\Es^{\Q^{*,1}}\left[\int_0^T {(\widetilde k^{\epsilon,i}_s)}^2
\widetilde \lambda^i_s ds \right]\le
\Es^{\Q^{*,1}}\left[{\left({1\over \bar K^\epsilon_T}
\right)}^2\right]\le \widetilde C_K
\end{split}
\end{equation}
From the uniqueness of the representation theorem and
(\ref{utildif}), we obtain,   for $1\le i\le d$, $-\Delta_\epsilon v^i_t=\ln\left[1+{\bar k^{\epsilon,i}_t\over\bar K^\epsilon_{t^-}} \right]$  and $\Delta_\epsilon v^i_t=\ln\left[1+\widetilde k^{\epsilon,i}_t K^\epsilon_{t^-}\right]$.Therefore we find $ \exp(|\Delta_\epsilon v^i_t|)-1\le {|\bar k^{\epsilon,i}_t|\over \bar K^{\epsilon}_{t^-}}+|\widetilde k^{\epsilon,i}_t|\bar K^{\epsilon}_{t^-}$. Moreover we have:
\begin{eqnarray*}
\Es^{\Q^{*,1}}\left[\int_0^T (e^{|\Delta_\epsilon v^i_s|}-1)\widetilde \lambda^i_s ds\right]&\le &\Es^{\Q^{*,1}}\left[ \int_0^T {|\bar k^{\epsilon,i}_s|\over \bar K^{\epsilon}_{s^-}} \widetilde \lambda^i_s ds +\int_0^T |\widetilde k^{\epsilon,i}_s|\bar K^{\epsilon}_{s^-} \widetilde \lambda^i_s ds \right]\\
&\le& \Es^{\Q^{*,1}}\left[\sup_{ t\le T}{1\over \bar K^{\epsilon}_t}\, \int_0^T|\bar k^{\epsilon,i}_s|\widetilde \lambda^i_s ds+ \sup_{t\le T} \bar K^{\epsilon}_t  \int_0^T |\widetilde k^{\epsilon,i}_s|\widetilde \lambda^i_s ds \right]\end{eqnarray*}
\noindent Using Cauchy-Schwartz inequality, we find that:
\begin{equation*}
\begin{split}
\Es^{\Q^{*,1}}\left[\int_0^T (e^{|\Delta_\epsilon v^i_s|}-1)\widetilde \lambda^i_s ds\right]&\le  {\left[\Es^{\Q^{*,1}}\left(\sup_{ t\le T}{1\over {(\bar K^{\epsilon}_t)}^2} \,\int_0^T \widetilde \lambda^i_s ds \right) \, \right]}^{1\over 2}{\left(\Es^{\Q^{*,1}}\int_0^T {|\widetilde k^{\epsilon,i}_t|}^2 \widetilde \lambda^i_s ds\right)}^{1\over 2}\\&+{\left[\Es^{\Q^{*,1}} \left(\sup_{t \le T}{(\bar K^{\epsilon}_t)}^2\,\int_0^T \widetilde \lambda^i_s ds\right)\right]}^{1\over 2} {\left(\Es^{\Q^{*,1}}\int_0^T {|\widetilde k^{\epsilon,i}_s|}^2 \widetilde \lambda^i_s ds\right)}^{1\over 2}.
\end{split}
\end{equation*}
We prove now that $\int_0^T \widetilde \lambda^i_t
dt$ is square integrable under the probability ${\Q^{*,1}}$. We write first the expression under $\P$ using Bayes's formula:
\begin{equation*}
\begin{split}
\Es^{\Q^{*,1}}{\left[\int_0^T \widetilde \lambda^i_t
dt\right]}^2=\Es^{\P}{\left[Z^{\Q^{*,1}}_T\int_0^T e^{-
v^{1,i}_s}\lambda^i_s ds \right]}^2 &\le \Es^\P\left[
{\big(Z^{\Q^{*,1}}_T\big)}^2{\int_0^T \lambda^i_s ds }{\int_0^T
e^{-2v^{1,i}_s}\lambda^i_s ds }\right]
\end{split}
\end{equation*}
\noindent then we use Cauchy-Schwartz inequality to find the following estimates:
\begin{equation*}
\begin{split}
 \Es^{\Q^{*,1}}{\left[\int_0^T \widetilde \lambda^i_t
dt\right]}^2 \le c_1{\left(\Es{\big(Z^{\Q^{*,1}}_T\big)}^4\right)}^{1\over 2}\left[\Es{\left(\int_0^T e^{-4 v^{1,i}_s}\lambda^i_s ds\right)}\right]^{1\over 2}
\end{split}
\end{equation*} where  we make use several times of Assumption A\ref{hyp2}-iii).
Moreover, we can see   that
\begin{equation*}
\begin{split}
\Es\left[\int_0^T e^{-4 v^{1,i}_s}\lambda^i_s ds\right]&=\Es\left[ \int_0^T {(e^{-v^{1,i}_s}-1+1)}^4\lambda^i_s ds\right] \le  16 \Es\left[ \int_0^T {(e^{-v^{1,i}_s}-1)}^4\lambda^i_s ds +\int_0^T \lambda^i_s ds \right].
\end{split}
\end{equation*}
 Therefore,  since the martingale $-M^{1,c}+\int_0^.\sum_{i=1}^d \left(e^{-v^{1,i}_t}-1\right)dN^i_t$ belongs to
 $ \mathbf{{L}}^p(\P)$, and by assumption A\ref{hyp2}-iii) again,  we conclude that
 $\Es^\P\left[ \int_0^T {(e^{-v^{1,i}_s}-1)}^p\lambda^i_s ds\right]<+\infty$ for any $p\ge 1$.
 Moreover since  $Z^{\Q^{*,1}}\in \mathbf{{L}}^p(\P)$, we get that
 $\Es^{\Q^{*,1}}\left[\int_0^T \widetilde\lambda^i_s ds \right] < \infty $.
\noindent then using again Cauchy inequality:
\begin{eqnarray*}
\Es^{\Q^{*,1}}\left[\int_0^T (e^{|\Delta_\epsilon
v^i_s|}-1)\widetilde \lambda^i_s ds\right]&\le &C
\Big(\Es^{\Q^{*,1}} \Big[\sup_{0\le t\le T}{1\over {(\bar
K^{\epsilon}_t)}^4}\Big]\Big)^{1\over 2}\Big(\Es^{\Q^{*,1}}\Big[\int_0^T
{|\bar k^{\epsilon,i}_s|}^2 \widetilde \lambda^i_s
ds\Big]\, \Big)^{1\over 2}\\&&+ C \, \Big(\Es^{\Q^{*,1}} \Big[\sup_{0\le t\le
T}{(\bar K^{\epsilon}_t)}^4 \Big] \Big)^{1\over 2}
\, \left(\Es^{\Q^{*,1}}\Big[\int_0^T {|\widetilde
k^{\epsilon,i}_s|}^2 \widetilde \lambda^i_s ds\Big]\right)^{1\over
2}
\end{eqnarray*}
From $(\ref{barK})$ and $(\ref{barsaut})$, we deduce that there
exists a constant $C_2>0$ such that:
$$\Es^{\Q^{*,1}}\left[\int_0^T (e^{|\Delta_\epsilon v^i_s|}-1)\widetilde \lambda^i_s ds\right]\le C_2$$
and then using the expansion of the functional $x\rightarrow e^x$ we get:
$$\Es^{\Q^{*,1}}\left[\int_0^T |\Delta_\epsilon v^i_s|^p\widetilde \lambda^i_s ds \right]\le C_2 p!.$$
 In order to  conclude the proof   of the lemma, it remains to
establish that there exists a constant $C_1$ satisfying:
$$\Es^{\Q^{*,1}}[\langle \Delta_\epsilon \widetilde M^{c}\rangle_T]\le C_1.$$
\noindent
\\[0.2cm]
\noindent \textit{First case}: $(c^2,\psi^2)\ge (c^1,\psi^1)$, then
$U(c^1+\epsilon(c^2-c^1))\ge U(c^1)$ and $\bar
U(\psi^1+\epsilon(\psi^2-\psi^1))\ge \bar U(\psi^1)$. From
Proposition \ref{estimation}, it  follows that:
\begin{equation*}
\begin{split}
&\Es^{\Q^{*,1}}\Big[ \sup_{0\le t\le T}|V^\epsilon_t-V^1_t|^2+
\langle \widetilde  M^{ \epsilon,c} -\widetilde M^{ 1,c} \rangle_T +\sum_{i=1}^d\int_0^T {(v^{\epsilon,i}_s-v^{1,i}_s)}^2\widetilde \lambda^i_s ds\Big]\\
&\le \Es^{\Q^{*,1}}\Big[{[\bar
U(\psi+\epsilon(\psi^2-\psi^1))-\bar U(\psi^1)]}^2+ \int_0^T
{[U(c^1_s+\epsilon(c^2_s-c^1_s))-U(c^1_s)]}^2ds\Big]
\end{split}
\end{equation*}
Since
$$ 0\le U(c^1_t+\epsilon(c^2_t-c^1_t))-U(c^1_t)\le \epsilon U'(c^1_t)(c^2_t-c^1_t)$$ and $$  0\le\bar U(\psi^1+\epsilon(\psi^2-\psi^1))-\bar U(\psi^1)\le\epsilon\bar U'(\psi^1)(\psi^2-\psi^1).$$
we get:
\begin{equation*}
\begin{split}
\Es^{\Q^{*,1}}\Big[\sup_{0\le t\le T}{|\Delta_\epsilon V_t|}^2+
\langle \Delta_\epsilon \widetilde M^{ c}\rangle_T +\int_0^T
{(\Delta_\epsilon v^i_s)}^2 \widetilde\lambda^{i}_sds \Big]&\le
\Es^{\Q^{*,1}}\Big[{(\bar
U'(\psi^1))}^2{(\psi^2-\psi^1)}^2 \\
& +\int_0^T {(U'(c^1_s))}^2{(c^2_s-c^1_s)}^2ds\Big] \\
\end{split}
\end{equation*}
 The process  $Z^{\Q^{*,1}}$ belongs to $\mathbf{{L}}^p(\P)$; moreover  $U'(\psi^1)(\psi^2-\psi^1) \in \mathbf{L}^{\exp}$ and $U'(c^1_s)(c^2_s-c^1_s)\in {D}^{\exp}_1$  since  $(c^1,\psi^1),(c^2,\psi^2)\in {\cal A}(x)$. It follows that  there exists a constant $C>0$
 such that:
\begin{equation*}
\begin{split}
\Es^{\Q^{*,1}}\Big[\sup_{0\le t\le T}{|\Delta_\epsilon V_t|}^2+
\langle \Delta_\epsilon \widetilde M^{c}\rangle_T +\sum_{i=1}^d
\int_0^T {(\Delta_\epsilon v^i_s)}^2 \widetilde\lambda^{i}_sds
\Big]&\le C
\end{split}
\end{equation*}\\\\ \noindent
\textit{ Second case}: $(c^2,\psi^2)\le (c^1,\psi^1)$. We first
prove   that for all $t\in[0,T]$, $\bar K^\epsilon_t\ge 1$. Let us
recall that:
 $$\bar K^\epsilon_t=\exp\left(-\Delta_\epsilon V_t+\int_0^t (\delta_s \Delta_\epsilon V_s-\Delta_\epsilon U_s)ds\right)$$
 \noindent Define the process $X$ as
$$X_t=-\Delta_\epsilon V_t +\int_0^t (\delta_s\Delta_\epsilon V_s -\Delta_\epsilon U_s)ds,\qquad 0\le t\le T.$$
\noindent From integration by part formula, we get:
\begin{equation*}
\begin{split}
S^\delta_tX_t&=-\Delta_\epsilon V_0 -\int_0^t S^\delta_s d\Delta_\epsilon V_s-\int_0^t S^\delta_s \Delta_\epsilon U_s ds\\
\end{split}
\end{equation*}
\noindent Since the process $\delta$ is positive and  bounded, there exists
a constant $L>0$ such that $S^\delta<L<1$.   It follows that:
$$S^\delta_t X_t\ge (-1+L)\Delta_\epsilon V_0 -L\Delta_\epsilon V_t-\int_0^t S^\delta_s \Delta_\epsilon U_s ds$$
 Note that,  for all $t\in [0,T]$,
$\Delta_\epsilon U_t\le 0$
since $(c^2,\psi^2)\le (c^1,\psi^1)$  and using comparison theorem  $\Delta_\epsilon V_t\le 0$.

Therefore, for all $t\in [0,T],X_t\ge 0$. Finally, $\bar
K^\epsilon_t\ge 1$.
\\  In the second step of the proof, we give the dynamics of the process $\bar K^\epsilon$ using It\^o's calculus:
$${d \bar K^\epsilon_t}=  \bar K^\epsilon_{t-}\left(-d \Delta_\epsilon \widetilde M^{c}_t +\sum_{i=1}^d (e^{-{(v^\epsilon_t-v^1_t)\over\epsilon}}-1)d\widetilde N^i_t +dA_t\right)$$
where $A$ is an increasing process. Since $\bar K^\epsilon$ is a
positive $\Q^{*,1}$-submartingale, we obtain from (\ref{barK}) and
$\bar K^\epsilon_t\ge 1$ :
$$\Es^{\Q^{*,1}}\left[\langle  \Delta_\epsilon \widetilde M^{c}\rangle_T\right]\le \Es^{\Q^{*,1}}\left[\int_0^T {(\bar K^\epsilon_t)}^2 d\langle \Delta_\epsilon \widetilde M^c\rangle _t\right]\le \Es^{\Q^{*,1}}\left[{(\bar K^\epsilon_T)}^2\right]\le C_K$$
then we conclude:
$$\Es^{\Q^{*,1}}\left[\langle  \Delta_\epsilon \widetilde M^{c}\rangle_T\right]\le C_K.$$
\noindent
 Finally, by using concavity property we have shown that: $|\Delta_\epsilon V_t|\le |V^2_t-V^1_t|$, $\hbox{ for all } t\in [0,T]$, then:
$$\Es^{\Q^{*,1}}\left[\sup_{t\in [0,T]}{|\Delta_\epsilon V_t|}^2\right]\le \Es^{\Q^{*,1}}\left[\sup_{t\in [0,T]}{|V^2_t-V^1_t|}^2\right]\le 2\, \Es^{\Q^{*,1}}\left[ \sup_{t\in [0,T]}{|V^1_t|}^2+\sup_{t\in [0,T]}{|V^2_t|}^2\right]$$
\noindent
Therefore, since the process $V^1,V^2\in D^{\exp}_0$ and  $Z^{\Q^{*,1}}$ belongs to $\bf{L}^p,$ we get by  using Cauchy Schwarz inequality  that there exists a constant $C$ such that:
$$\Es^{\Q^{*,1}}\left[\sup_{t\in [0,T]}{|\Delta_\epsilon V_t|}^2\right]\le C.$$

$\hfill\Box$

\paragraph{Proof of Theorem \ref{difdym}:}
Let recall first the equality :
$${1\over 2  }(\langle M^{ \epsilon,c }\rangle-\langle M^{ 1,c }\rangle)={1\over 2}\langle M^{ \epsilon,c}-
M^{ 1,c}\rangle +\langle M^{ \epsilon,c},M^{ 1,c}\rangle-\langle
M^{1,c}\rangle,$$ then the equation \eqref{BSDE:approximation} may
be written as:
\begin{equation*}
\begin{split}
&\Delta_\epsilon V_t-\int_0^t(\delta_s \Delta_\epsilon
V_s-\Delta_\epsilon U_s)ds  = {1\over \epsilon}\Big({1\over
2}\langle M^{\epsilon,c}-M^{1,c}\rangle_t +\langle
M^{\epsilon,c},M^{1,c}\rangle_t-\langle M^{1,c}\rangle_t
\Big)\\
& +\sum_{i=1}^d \int_0^t \big[{1\over
\epsilon}(e^{-v^{\epsilon,i}_s}-e^{-v^{1,i}_s})+e^{-v^{1,i}}\Delta_\epsilon
v^i_s \big]\, \lambda^i_s  ds +\Delta_\epsilon M^{c}_t
 +\sum_{i=1}^d\int_0^t \Delta_\epsilon v^i_s
(dN^i_s-(e^{v^{1,i}}-1)\lambda^i_s ds)\\
&={1\over 2
\epsilon}\langle M^{\epsilon,c}-M^{1,c}\rangle_t +\sum_{i=1}^d\int_0^t
  \Big[{1\over
\epsilon}(e^{-v^{\epsilon,i}_s}-e^{-v^{1,i}_s})+e^{-v^{1,i}}\Delta_\epsilon
v^i_s \Big]\lambda^i_s ds\\
&+(\Delta_\epsilon M^{c}_t +\langle
\Delta_\epsilon M^{c},M^{1,c}\rangle_t)+\sum_{i=1}^d\int_0^t
\Delta_\epsilon v^i_s (dN^i_s-(e^{-v^{1,i}}-1)\lambda^i_s ds).
\end{split}
\end{equation*}
\\
\noindent
By Girsanov's theorem, the processes $ \Delta_\epsilon \widetilde M^{c}:=\Delta_\epsilon
M^{c} +\langle \Delta_\epsilon M^{c},M^{1,c}\rangle$ and
$\widetilde N^i:=N^i-\int_0^.(e^{-v^{1,i}_s}-1)\lambda^i_s ds$ are
$\Q^{1,*}-$ martingales. It follows that the process:
$${(\Delta_\epsilon V_t-\int_0^t(\delta_s \Delta_\epsilon
V_s-\Delta_\epsilon U_s)ds)}_{t\ge 0}$$ \noindent is a
$\Q^{1,*}$-submartingale and we have the following decomposition:
\begin{equation}\label{equation:Q-approximation}
\begin{split}
&\Delta_\epsilon V_t-\int_0^t(\delta_s \Delta_\epsilon
V_s-\Delta_\epsilon U_s)ds={\epsilon\over 2} \langle
\Delta_\epsilon \widetilde M^{c}\rangle_t+ \sum_{i=1}^d\int_0^t
\big[{1\over
\epsilon}(e^{-v^{\epsilon,i}_s}-e^{-v^{1,i}_s})\\&+e^{-v^{1,i}}\Delta_\epsilon
v^i_s \big]\, \lambda^i_s  ds +\Delta_\epsilon \widetilde
M^{c}_t +\sum_{i=1}^d\int_0^t  \Delta_\epsilon v^i_s d\widetilde
N^i_s.
\end{split}
\end{equation}
\noindent Moreover, using the uniform estimate \eqref{priori:estimate}, we get:
\begin{equation}
\label{M:approximation}
\lim_{\epsilon\rightarrow 0}  \Es^{\Q^{*,1}}\left({\epsilon\over 2} \langle \Delta_\epsilon \widetilde M^{c}\rangle_T\right)\le C_p\lim_{\epsilon\rightarrow 0} {\epsilon\over 2}=0,
\end{equation}
and using the expansion of the functional $x\rightarrow e^x$, we get:
\begin{eqnarray*}
0&\le& \lim_{\epsilon \rightarrow 0}
\Es^{\Q^{*,1}}\left(\int_0^T[ \frac{e^{-v^{\epsilon,i}_s}-e^{-v^{1,i}_s}}{\epsilon}+e^{-v^{1,i}_s}\Delta_\epsilon
v^i_s]\lambda^i_s ds\right)=\lim_{\epsilon\rightarrow
0}\Es^{\Q^{*,1}}\left(\int_0^T \sum_{p=2}^{+\infty}{\epsilon^{p-1}\over
p!}{(\Delta_\epsilon v^i_s)}^p\widetilde\lambda^{i}_s ds
\right)\\&& \le \sum_{p=2}^{\infty}{\epsilon^{p-1}} \Es^{\Q^{1,*}}\left( \int_0^T {|\Delta_\epsilon v^i_s|^p\over p!}\widetilde\lambda^{i}_s ds\right)
 \le \sum_{p=2}^{\infty}C {\epsilon^{p-1}}={C \epsilon \over 1-\epsilon},
 \end{eqnarray*}
thus, passing to the limit as $ \epsilon \to 0$, we conclude  that:
\begin{equation}
\label{v:approximation}
\lim_{\epsilon \rightarrow 0} \Es^{Q^{*,1}}\left(\int_0^T[{1\over \epsilon}(e^{-v^{\epsilon,i}_s}-e^{-v^{1,i}_s})+e^{-v^{1,i}}\Delta_\epsilon v^i_s]\lambda^i_s ds\right)=0,\qquad 1\le i\le d.
\end{equation}
\noindent Moreover, the estimate \eqref{priori:estimate} ensures that the
sequence $ (\Delta_\epsilon V, \Delta_\epsilon \widetilde{M^{ c}},
\Delta_\epsilon v)_{\epsilon >0}$ is bounded in $
\mathcal{H}^2(\P) \times \mathcal{M}_0^2(\P) \times
 \mathcal{L}^2(\lambda, \P) $.  As a consequence,  we can
extract a subsequence $ (\Delta_{\epsilon_k} V,
\Delta_{\epsilon_k} M^{ c}, \Delta_{\epsilon_k} v )_{k \in
\mathbb{N}}$ which converges weakly in $\mathcal{H}^2(\P) \times
\mathcal{M}_0^2(\P) \times  \mathcal{L}^2(\lambda, \P) $ and by
Banach-Mazur Lemma, one may construct a sequence $
(\widehat{V}^{\epsilon}, \widehat{M}^{\epsilon, c},
\widehat{v}^{\epsilon})_{\epsilon >0}$ of convex combinations of
elements in  $(\Delta_{\epsilon_k} V, \Delta_{\epsilon_k}
\widetilde{M^{ c}}, \Delta_{\epsilon_k} v )_{k \in N}$ of the form
$$ \widehat{V}^\epsilon :=\sum_{j=1}^{N_{\epsilon} } \alpha^{\epsilon_j} \Delta_{\epsilon_j} V, \quad
 \widehat{M}^{\epsilon,c} :=\sum_{j=1}^{N_{\epsilon} } \alpha^{\epsilon_j} \Delta_{\epsilon_j} \widetilde{M^c}, \quad
  \widehat{v}^\epsilon :=\sum_{j=1}^{N_{\epsilon} } \alpha^{\epsilon_j} \Delta_{\epsilon_j} v
 $$
 such that
 $ (\widehat{V}^{\epsilon}, \widehat{M}^{\epsilon,c}, \widehat{v}^{\epsilon})_{\epsilon >0}$
 converges strongly in $\mathcal{H}^2(\P) \times \mathcal{M}_0^2(\P) \times
 \mathcal{L}^2(\lambda, \P) $ to  $(\partial_\epsilon
V,\partial_\epsilon \widetilde M^{c},\partial_\epsilon v)$.
Moreover, the triple $ (\widehat{V}^{\epsilon},
\widehat{M}^{\epsilon, c}, \widehat{v}^{\epsilon})$ satisfies the
BSDE \eqref{equation:Q-approximation} associated with  $ (
\widehat{U}^\epsilon, \widehat{\bar U}^\epsilon)$ where $$
\widehat{U}^\epsilon := \sum_{j=1}^{N_{\epsilon} }
\alpha^{\epsilon_j} \Delta_{\epsilon_j} U, \quad \widehat{\bar
U^\epsilon} := \sum_{j=1}^{N_{\epsilon} } \alpha^{\epsilon_j}
\Delta_{\epsilon_j} \bar U .$$ Therefore,  passing  to the limit
in this equation, thanks to   \eqref{M:approximation},
\eqref{v:approximation} and the dominated convergence theorem, we
get that  $(\partial_\epsilon V,\partial_\epsilon \widetilde
M^{c},\partial_\epsilon v)$ solves the BSDE
$$d\partial_\epsilon V_t=(\delta_t \partial_\epsilon
V_t-U'(c^1_t)(c^2_t-c^1_t))dt +d\partial_\epsilon \widetilde
M^{c}_t+\sum_{i=1}^d \partial_\epsilon v^i_t d\widetilde N^i_t,\quad \partial_\epsilon  V_T=\bar U'(\psi^1)(\psi^2-\psi^2).
$$
\noindent Therefore ${(S^\delta_t \partial_\epsilon V_t+\int_0^t S^\delta_s
U'(c^1_s)(c^2_s-c^1_s)ds)}_{t\ge 0}$ is a $\Q^{*,1}$ martingale
which can  be written as:
$$ S^\delta_t\partial_\epsilon V_t+\int_0^t S^\delta_s U'(c^1_s)(c^2_s-c^1_s) ds=\Es^{\Q^{1,*}}\Big[S^\delta_T\partial_{\epsilon} V_T+\int_0^T S^\delta_s U'(c^1_s)(c^2_s-c^1_s)ds \Big\vert {\cal G}_t\Big].$$
Hence we get:
$$\partial_\epsilon V_t=\Es^{\Q^{1,*}}\Big[{S^\delta_T\over S^\delta_t}\bar U'(\psi^1)(\psi^2-\psi^1)+\int_t^T {S^\delta_s\over S^\delta_s} U'(c^1_s)(c^2_s-c^1_s)ds \Big\vert {\cal G}_t\Big].$$
$\hfill\Box$

\newpage
\noindent 
Monique JEANBLANC, \\
Laboratoire de Math\'ematiques et Mod\'elisation d'\'Evry (LaMME),   Universit\'e d'\'Evry-Val-d'Essonne, UMR CNRS 8071\\
23, Boulevard de France,\\
 F-91037 Evry Cedex, FRANCE\\
 e-mail: monique.jeanblanc{\char'100}univ-evry.fr \\[0.3cm]
 Anis MATOUSSI,\\
   Université du Maine,\\
 Laboratoire Manceau de Mathématiques,\\
 Avenue Olivier Messiaen,
 F-72085 Le Mans Cedex 9, France \\
 and  \\
CMAP,  Ecole Polytechnique, Palaiseau \\
email : anis.matoussi@univ-lemans.fr
\\[0.3cm]
Armand NGOUPEYOU,\\
Banque des Etats de l'Afrique Centrale\\
Agence de Limbé, Service des Etudes\\
BP 15,\\
Limbé, Cameroun\\
ngoupeyou@beac.int\\

\end{document}